\input amstex 
\documentstyle{amsppt} 
\loadbold

\let\bs\boldsymbol
\magnification=1200
\hsize=5.75truein
\vsize=8.75truein 
\hcorrection{.25truein}
\loadeusm 
\let\scr\eusm
\loadeurm 


\define\aC{\text{\it a-}\scr C} 
\define\cind{\text{\it c-\/\rm Ind}} 

\define\Aut#1#2{\text{\rm Aut}_{#1}(#2)}
\define\End#1#2{\text{\rm End}_{#1}(#2)}
\define\Hom#1#2#3{\text{\rm Hom}_{#1}({#2},{#3})} 
\define\GL#1#2{\roman{GL}_{#1}(#2)}
\define\M#1#2{\roman M_{#1}(#2)}
\define\Gal#1#2{\text{\rm Gal\hskip.5pt}(#1/#2)}
\define\awr#1{\scr A^{\text{\rm m-wr}}(#1)} 
\define\asq#1{\scr A^{\square}(#1)} 
\define\Awr#1#2{\scr A^{\text{\rm m-wr}}_{#1}(#2)} 
\define\upr#1#2{{}^{#1\!}{#2}} 
\define\pre#1#2{{}_{#1\!}{#2}} 
\define\N#1#2{\text{\rm N}_{#1/#2}} 
\define\Mid{\,\big|\,} 
\define\Swr#1{\scr S^\roman{wr}(#1)} 
\define\r#1{\text{\rm #1}} 
\let\ge\geqslant
\let\le\leqslant 
\let\wt\widetilde 
\let\ups\upsilon

\let\vf\varphi
\let\vF\varPhi 
\let\vG\varGamma 
 
\let\vL\varLambda 
\let\vO\varOmega 
\let\vp\varpi 
 
\let\vP\varPi

\let\vt\vartheta

\document \baselineskip=14pt \parskip=4pt plus 1pt minus 1pt
\topmatter \nologo \nopagenumbers
\title\nofrills 
Explicit local Jacquet-Langlands correspondence: \\ the non-dyadic wild case
\endtitle 
\rightheadtext{Wild correspondence} 
\author 
Colin J. Bushnell and Guy Henniart 
\endauthor 
\leftheadtext{C.J. Bushnell and G. Henniart}
\affil 
King's College London and Universit\'e de Paris-Sud 
\endaffil 
\address 
King's College London, Department of Mathematics, Strand, London WC2R 2LS, UK. 
\endaddress
\email 
colin.bushnell\@kcl.ac.uk 
\endemail
\address 
Laboratoire de Math\'ematiques d'Orsay, Univ\. Paris-Sud, CNRS, Universit\'e
Paris-Saclay, 91405 Orsay, France.
\endaddress 
\email 
Guy.Henniart\@math.u-psud.fr 
\endemail 
\date September 2017 \enddate 
\abstract 
Let $F$ be a non-Archimedean locally compact field of residual characteristic $p$ with $p\neq 2$. Let $n$ be a power of $p$ and let $G$ be an inner form of the general linear group $\text{\rm GL}_n(F)$. We give a transparent parametrization of the irreducible, totally ramified, cuspidal representations of $G$ of parametric degree $n$. We show that the parametrization is respected by the Jacquet-Langlands correspondence, relative to any other inner form. This expresses the Jacquet-Langlands correspondence for such representations within a single, compact formula. 
\endabstract 
\endtopmatter 
\document \baselineskip=14pt \parskip=4pt plus 1pt minus 1pt 
\subhead 
1 
\endsubhead 
Let $F$ be a non-Archimedean local field of residual characteristic $p$. Let $n\ge1$ and let $G$ be an inner form of the general linear group $\GL nF$. In other words, there is a central simple $F$-algebra $A$, of dimension $n^2$, such that $G = A^\times$. Let $\asq G$ be the set of equivalence classes of {\it essentially square-integrable,\/} smooth, complex representations of $G$. Let $G'$ be another inner form of $\GL nF$. We study the canonical bijection 
$$ 
T^{G'}_G: \asq G @>{\ \ \approx\ \ }>> \asq{G'} 
$$ 
provided by the Jacquet-Langlands correspondence \cite{10}, \cite{1}. We make a narrow, but significant, contribution to the analysis of the correspondence in explicit terms. 
\par 
Let $\pi\in \asq G$ and let $d(\pi)$ be the {\it parametric degree\/} of $\pi$, in the sense of \cite{6}. Thus $d(\pi)$ is a positive integer dividing $n$. If $d(\pi) = n$, then $\pi$ is cuspidal. The converse holds if $G$ is the split group $\GL nF$ but not in general: for example, if $G = \GL1D$, where $D$ is a central $F$-division algebra of dimension $n^2$, any irreducible smooth representation $\pi$ of $G$ is cuspidal while $d(\pi)$ is arbitrary. The parametric degree is preserved by the Jacquet-Langlands correspondence. 
\par 
In this paper, we concentrate on the case where $\pi$ is of parametric degree $n$ and {\it totally wildly ramified.} This means that $n$ is a power of $p$ and, if $\chi\neq 1$ is an unramified character of $F^\times$, the twist $\chi\pi$ of $\pi$ is not equivalent to $\pi$. When $p\neq 2$, such representations admit a particularly transparent description (2.4 Proposition) that we can use to describe the correspondence via a compact explicit formula (6.2 Theorem). The case $p=2$ has sufficiently many distinctive features to merit a separate treatment that we defer for the time being. 
\par 
With this result to hand, the way is open to follow the framework of \cite{6} and \cite{8} (but without complications arising from the transfer factors of automorphic induction \cite{13}) to an explicit description of the Jacquet-Langlands correspondence for representations $\pi\in \asq G$ with $d(\pi) = n$. The more difficult case is that where $d(\pi) < n$ while $\pi$ is cuspidal. With the newly available Endo-class Transfer Theorem of \cite{22} and \cite{11} recalled below, that general case is substantially less mysterious than hitherto. However, it seems unlikely that one will resolve the question finally without the detail of the complementary special case treated here. The one fully known case of \cite{24}, \cite{7} indicates the prospect of some intriguing further subtlety. 
\subhead 
2 
\endsubhead 
For the convenience of the reader, we interpolate an outline summary of recent developments in the broader context. Again, $G$ is an inner form of $\GL nF$, but we impose no restriction on $n$ or $p$ for the time being. 
\par 
The papers \cite{3,17--21} of S\'echerre and collaborators contain a complete description of the representations $\pi \in\asq G$ in terms of {\it simple characters\/} and {\it simple types.} Allowing for a few novel features and a higher level of technical intricacy, it is parallel to the split case $G = \GL nF$ of \cite{9}. In particular, any $\pi\in \asq G$ contains a simple type and hence a simple character $\theta_\pi$. The representations $\pi$ that contain a given simple type are classified via a scheme following that of the split case \cite{19}. 
\par 
Simple characters, as a class, have a fundamental naturality property. Work\-ing at first in the split case of \cite{9}, let $\theta$ be a simple character in $G = \GL nF$. Thus $\theta$ is attached to a hereditary $\frak o_F$-order in the matrix algebra $\M nF$. If $\frak a'$ is a hereditary order in $\M{n'}F$ then, subject to minor combinatorial constraints, one can construct from $\theta$ a simple character $\theta'$ in $\GL{n'}F$, attached to $\frak a'$. We refer to $\theta'$ as a ``transfer'' of $\theta$. If, for $i=1,2$, we are given a simple character $\theta_i$ in $G_i = \GL{n_i}F$, one can always find an integer $n_3$ and a hereditary order $\frak a_3$ in $A_3 = \M{n_3}F$ that admits a transfer $\theta_i'$ of $\theta_i$, $i=1,2$. One knows from \cite{4} that if, for some choice of datum $(A_3,\frak a_3)$, the transfers $\theta_i'$ intertwine (and so are conjugate) in $A_3^\times$, then the same is true for all such choices. When this holds, one says that $\theta_1$ is {\it endo-equivalent to\/} $\theta_2$. Endo-equivalence is an equivalence relation on the class of all simple characters in all groups $\GL nF$, $n\ge1$. The set of endo-equivalence classes (endo-classes for short) is an arithmetic object of considerable interest: see section 6 of \cite{8} for an overview. 
\par 
The achievement of \cite{3} is an extension of this relation to the class of all simple characters in all inner forms of all $\GL nF$. Every endo-equivalence class, in this extended context, contains a simple character in some split group $\GL nF$, $n\ge1$. On such characters, the two notions of endo-equivalence are the same. For general $G$ and $\pi\in \asq G$, the endo-class of a simple character $\theta_\pi$ contained in $\pi$ is uniquely determined by $\pi$. One cannot avoid asking how this fundamental invariant behaves with respect to the Jacquet-Langlands correspondence. 
\proclaim{Endo-class Transfer Theorem \cite{22}, \cite{11}} 
For $i=1,2$, let $G_i$ be an inner form of $\GL nF$. Let $\pi_i\in \asq{G_i}$ and let $\theta_i$ be a simple character contained in $\pi_i$.  If $\pi_2 = T^{G_2}_{G_1}(\pi_1)$, then $\theta_1$ is endo-equivalent to $\theta_2$. 
\endproclaim 
The proof of this result takes an unexpected form. Let $\ell$ be a prime number different from $p$. In a series of papers including \cite{14} and \cite{15}, M\'\i nguez and S\'echerre develop a theory of $\ell$-modular representations of the inner forms $G$ of $\GL nF$ and of reduction, modulo $\ell$, of representations in characteristic zero. In \cite{16}, they show that reduction modulo $\ell$ is compatible with the Jacquet-Langlands correspondence. Using these results, for varying $\ell$, S\'echerre and Stevens show that the Endo-class Transfer Theorem holds in general {\it provided it holds when $d(\pi_i) = n$ and $\pi_1$ is totally ramified\/} \cite{22}, that is, $\chi\pi_1 \not\cong \pi_1$ when $\chi$ is a non-trivial unramified  character of $F^\times$. Using a neat device combining properties of certain simple characters, relative to unramified base field extension, and trace comparisons of a sort familiar from \cite{6} or \cite{8}, Dotto \cite{11} reduces to the split groups and so despatches the outstanding special case. That method also yields the relation between the simple  types  contained in corresponding representations $\pi_i$ of parametric degree $n$. 
\subhead 
3 
\endsubhead 
We return to the theme of this paper. From now on, $G$ is an inner form of $\GL nF$, where $n = p^r$, $r\ge1$. Beyond the very first stages, we assume $p\neq2$. A representation $\pi\in \asq G$ is {\it totally wildly ramified\/} if $\pi\not\cong \chi\pi$ for any unramified character $\chi\neq1$ of $F^\times$. Let $\awr G$ be the set of totally wildly ramified representations $\pi\in \asq G$ such that $d(\pi) = n$. A representation $\pi \in \awr G$ is of the form $\cind_{\r J}^G\,\vL$, where $(\r J,\vL)$ is an {\it extended maximal simple type\/} in $G$ \cite{20}. In general, the representation $\vL$ has dimension $p^s$, for an integer $s\ge1$ and its character (that is, trace) is not conveniently accessible. The strategy is to replace $(\r J,\vL)$ by a pair $(\r I,\lambda)$, in which $\r I$ is a certain canonical subgroup of $\r J$ and $\lambda$ is a character of $\r I$, extending the simple character in $\pi$: the constructions are given in sections 1 and 2. The representation $\pi$ is induced by any of these characters $\lambda$ that it contains. However, one may choose $\lambda$ to have a particular, explicit form: this ``standard form'' is written out in 2.4. It leads to an exact parametrization of those $\pi$, containing a fixed simple character $\theta$, in terms of characters of a field associated with $\theta$ (2.4 Proposition). 
\par 
The one-dimensional parameters $(\r I,\lambda)$ behave transparently with respect to finite, unramified base field extension and can be transferred, via such an extension, to an inner form $G'$ of $G$. This transfer process, set out in section 3, specializes to the standard transfer of simple characters as in \cite{3,17}, but it also suggests an explicit {\it parametric transfer\/} of representations between $\awr G$ and $\awr{G'}$, say $\pi\mapsto \pi'$. At this stage, the parametric transfer is not well-defined, apparently depending on a fairly random choice. However, choosing correctly, it does preserve the standard form of 2.4 (by 3.5 Proposition). 
\par 
There is a second approach (in section 4) obtained by passing to the completion of the maximal unramified extension $\wt F$ of $F$. The transfer process applies equally over $\wt F$ but, because of completeness, it can be achieved via a conjugation in the group of $\wt F$-points of $G$. This version of the parametric transfer is equally ill-defined but is equivalent to the first one (4.5 Proposition). However, the fact that it is given by a conjugation over the complete field $\wt F$ enables comparison of the characters of $\pi$ and $\pi'$. In section 5, we show that the characters $\roman{tr}\,\pi$, $\roman{tr}\,\pi'$ agree at enough  elements to ensure they are equal. The conclusion in section 6 is that {\it the parametric transfer, in either version, is the Jacquet-Langlands correspondence.} The Jacquet-Langlands correspondence therefore respects the standard form of 2.4, whence it is expressed as a compact and transparent formula (6.2 Theorem). 
\par 
This paper up-dates and supersedes the relevant parts of our earlier work \cite{5}. That concerned only the relation between $\GL nF$ and $\GL1D$, where $D$ is a central $F$-division algebra of dimension $n^2$. Not enough of the general machinery of \cite{3,17--21} was available at that time, so \cite{5} could only rely on \cite{2}. However, a lot of the effort in \cite{5} is centred on $\GL nF$, and is used to ease our task here. On the other hand, the step from $\GL1D$ to a general inner form requires some effort and re-organization of the detail into a more efficient and flexible form. 
\head 
1. The Lagrangian subgroup
\endhead 
Let $A$ be a central simple $F$-algebra of dimension $n^2$, $n = p^r$, and set $G = A^\times$. Let $\awr G$ denote the set of equivalence classes of irreducible, smooth, complex representations of $G$ that are cuspidal, totally ramified and of parametric degree $n$. In this and the following section, we describe the elements $\pi$ of $\awr G$ as representations induced from a canonical family of {\it characters\/} of open, compact modulo centre, subgroups of $G$. 
\par 
We recall something of the {\it simple characters in\/} $A$, following a simplified version of the foundational account of \cite{17}: since we deal only with a very special case, the more elaborate technical structures of \cite{17} are not needed here. From there, we develop a modified version of the method of \cite{5}. In this section, we allow the possibility $p=2$. 
\subhead 
1.1 
\endsubhead 
Let $\frak a$ be a {\it minimal hereditary $\frak o_F$-order\/} in $A$. Thus $\frak a$ is a principal order: if $\frak p$ is the Jacobson radical of $\frak a$, there exists $\vP \in G$ such that $\frak p = \vP\frak a = \frak a\vP$. Any two minimal hereditary orders in $A$ are $G$-conjugate. 
\par 
We use the concept of {\it simple stratum in\/} $A$, following \cite{17}. 
\definition{Notation} 
Let $\Swr{\frak a}$ be the set of elements $\beta$ of $G$ satisfying the following conditions. 
\roster 
\item There is an integer $l>0$ such that the quadruple $[\frak a,l,0,\beta]$ is a simple stratum in $A$. 
\item The field extension $F[\beta]/F$ is of degree $n$. 
\endroster 
\enddefinition 
Since $\frak a$ is minimal and $n$ is a power of $p$, these conditions imply that $F[\beta]/F$ is {\it totally wildly ramified.} The order $\frak a$ is stable under conjugation by $F[\beta]^\times$ --- one says that $\frak a$ is {\it $F[\beta]$-pure} --- and $\frak a$ is the unique hereditary order in $A$ with this property. The integer $l$ is given by $\beta^{-1}\frak a = \frak p^l$. 
\proclaim{Proposition} 
Let $\beta\in \Swr{\frak a}$ and write $E = F[\beta]$. Let $B$ be a central simple $F$-algebra of dimension $n^2$, let $\iota:E\to B$ be an $F$-embedding, and let $\frak b$ be an $\iota E$-pure hereditary $\frak o_F$-order in $B$. The stratum $[\frak b,l,0,\iota\beta]$ is then simple, the order $\frak b$ is minimal and $\iota\beta\in \Swr{\frak b}$. 
\endproclaim 
\demo{Proof} 
See Proposition 2.25 of \cite{17}. \qed 
\enddemo 
\subhead 
1.2 
\endsubhead 
Let $\beta\in \Swr{\frak a}$. Following \cite{17}, the simple stratum defined by $\frak a$ and $\beta$ gives rise to a pair of $\frak o_F$-orders in $A$, 
$$ 
\frak H(\beta,\frak a) \i \frak J(\beta,\frak a) \i \frak a, 
$$ 
and families of open subgroups 
$$ 
\aligned 
H^k(\beta,\frak a) &= 1+\frak H(\beta,\frak a)\cap \frak p^k, \\ 
J^k(\beta,\frak a) &= 1+\frak J(\beta,\frak a)\cap \frak p^k, 
\endaligned \qquad k\ge1, 
$$ 
of the principal unit group $U^1_\frak a = 1{+}\frak p$. 
\par 
Fix a character $\psi^F$ of $F$ that is trivial on $\frak p_F$ but not trivial on $\frak o_F$: one says that $\psi^F$ {\it is of level one.} As in \cite{17}, use $\psi^F$ to define the set $\scr C(\frak a,\beta,\psi^F)$ of {\it simple characters\/} of $H^1(\beta,\frak a)$. 
\par 
We write $\psi^A = \psi^F\circ \roman{tr}_A$, where $\roman{tr}_A:A\to F$ is the reduced trace map. For $\alpha\in A$, we define a function $\psi^A_\alpha$ by 
$$ 
\psi^A_\alpha(x) = \psi^A(\alpha(x{-}1)), \quad x\in A. 
\tag 1.2.1 
$$ 
\subhead 
1.3 
\endsubhead 
We recall from \cite{17} {\it passim\/} the behaviour of these structures relative to unramified base field extension. 
\par 
Let $K/F$ be a finite, {\it unramified\/} field extension. The $K$-algebra $A_K = A\otimes_F K$ is central simple of dimension $n^2$. Set $G_K = A_K^\times$. The ring $\frak a_K = \frak a\otimes_{\frak o_F}\frak o_K$ is a minimal hereditary $\frak o_K$-order in $A_K$, with Jacobson radical $\frak p_K = \frak p\otimes_{\frak o_F}\frak o_K$. We habitually identify $A$ with the subring $A\otimes 1$ of $A_K$. 
\proclaim{Proposition} 
Let $K/F$ be a finite unramified field extension and let $\beta\in \Swr{\frak a}$. 
\roster 
\item 
 The element $\beta\otimes1$ of $G_K$ lies in $\scr S^\roman{wr}(\frak a_K)$ and 
$$ 
\align 
\frak H(\beta\otimes 1,\frak a_K) &= \frak H(\beta,\frak a) \otimes_{\frak o_F}\frak o_K, \\ 
H^k(\beta,\frak a) &= H^k(\beta\otimes 1,\frak a_K)\cap G, \quad k\ge1. 
\endalign 
$$ 
Similarly for the $J$-groups. 
\item 
Let $\psi^K$ be a character of $K$, of level one, such that $\psi^K\Mid F = \psi^F$. If $\theta\in \scr C(\frak a_K,\beta\otimes 1,\psi^K)$, the character $\theta^F = \theta\Mid H^1(\beta,\frak a)$ lies in $\scr C(\frak a,\beta,\psi^F)$. The restriction map 
$$ 
\align 
\scr C(\frak a_K,\beta\otimes 1,\psi^K) &\longrightarrow \scr C(\frak a,\beta,\psi^F), \\ 
\theta &\longmapsto \theta^F, 
\endalign 
$$ 
is surjective. 
\endroster 
\endproclaim 
\demo{Proof} 
If the degree $[K{:}F]$ is divisible by $n$, then $A_K \cong \M nF$ and all assertions follow directly from the definitions in \cite{17}, particularly 3.3. The general case then follows by transitivity. \qed 
\enddemo 
From now on, we follow convention and write $\beta = \beta\otimes 1 \in A_K$. 
\subhead 
1.4 
\endsubhead 
Denote by $\r K_\frak a$ the group of $g\in G$ for which $g\frak ag^{-1} = \frak a$. Equivalently, $\r K_\frak a$ is the $G$-normalizer of $U_\frak a = \frak a^\times$. It is generated by $U_\frak a$ and any element $\vP$ such that $\vP\frak a = \frak p$. 
\par 
Let $\beta\in \Swr{\frak a}$ and let $\theta\in \scr C(\frak a,\beta,\psi^F)$. The $G$-normalizer of $\theta$ is the group 
$$ 
\r J(\beta,\frak a) = F[\beta]^\times J^1(\beta,\frak a). 
$$ 
In particular, $\r J(\beta,\frak a)$ is an open subgroup of $\r K_\frak a$ that does not depend on the choice of $\theta \in \scr C(\frak a,\beta,\psi^F)$. An element $g$ of $G$ intertwines the character $\theta$ if and only if $g\in \r J(\beta,\frak a)$. (For these facts, see \cite{17} Th\'eor\`eme 3.50.) 
\par 
The point of the section is to construct a canonical open subgroup $I^1(\beta,\frak a)$ of $G$, lying between $J^1(\beta,\frak a)$ and $H^1(\beta,\frak a)$. The group 
$$ 
\scr J^1(\beta,\frak a) = J^1(\beta,\frak a)/H^1(\beta,\frak a) \cong \frak J^1(\beta,\frak a)/\frak H^1(\beta,\frak a) 
$$ 
is a vector space over the finite residue field $\Bbbk_F = \Bbbk_{F[\beta]}$. In particular, it is a vector space over the field $\Bbb F_p$ of $p$ elements. Let $\theta\in \scr C(\frak a,\beta,\psi^F)$. Using the commutator convention $[x,y] = x^{-1}y^{-1}xy$,  the pairing 
$$ 
(x,y) \longmapsto \theta([x,y]), \quad x,y\in J^1(\beta,\frak a), 
\tag 1.4.1 
$$ 
induces an $\Bbb F_p$-bilinear form on $\scr J^1(\beta,\frak a)$. This form is nondegenerate and alternating \cite{17} Th\'eor\`eme 3.52. \proclaim{Lemma} 
The pairing \rom{(1.4.1)} satisfies 
$$ 
\theta([1{+}x,1{+}y]) = \psi^A_\beta(1{-}xy{+}yx), 
$$ 
for $x,y \in \frak J^1(\beta,\frak a)$ and $\theta\in \scr C(\frak a,\beta,\psi^F)$. 
\endproclaim 
\demo{Proof} 
When $A\cong \M nF$, the result is 6.1 Proposition of \cite{5}. We reduce the general case to that one. Let $K/F$ be a finite unramified extension such that $A_K \cong \M nK$. By 1.3 Proposition, there exists $\theta_K\in \scr C(\frak a_K,\beta,\psi^K)$ such that $\theta = \theta_K\Mid H^1(\beta,\frak a)$. For $x,y\in \frak J^1(\beta,\frak a) \i \frak J^1(\beta,\frak a_K)$, we have 
$$ 
\theta([1{+}x,1{+}y]) = \theta_K([1{+}x,1{+}y]) = \psi^{A_K}_\beta(1{-}xy{+}yx) 
$$ 
{\it loc\. cit.} On the other hand, $\psi^A_\beta = \psi^{A_K}_\beta\Mid A$ so 
$$ 
\theta([1{+}x,1{+}y]) =  \psi^A_\beta(1{-}xy{+}yx) , 
$$ 
as required. \qed 
\enddemo 
The pairing (1.4.1) is thus independent of the choice of $\theta \in \scr C(\frak a,\beta,\psi^F)$: we name it $h_\beta$, or $h_\beta^F$ when we need to specify the base field. 
\par 
Let $k\ge1$ be an integer. Define $\scr J^k = \scr J^k(\beta,\frak a)$ as the image of $J^k(\beta,\frak a)$ in $\scr J^1(\beta,\frak a)$. By Th\'eor\`eme 3.52 of \cite{17}, the pairing $h_\beta$ is nondegenerate on $\scr J^k(\beta,\frak a)$, $k\ge1$. So, for each $k\ge1$, there is a unique subspace $\scr U^k(\beta,\frak a)$ of $\scr J^1(\beta,\frak a)$ such that 
$$ 
\scr J^k(\beta,\frak a) = \scr U^k(\beta,\frak a) \perp \scr J^{k+1}(\beta,\frak a), 
$$ 
the sum being orthogonal with respect to the alternating form $h_\beta$. One has $\scr J^k(\beta,\frak a) \neq \scr J^{k+1}(\beta,\frak a)$ if and only if $2k$ is a jump of the stratum $[\frak a,l,0,\beta]$, so $\scr U^k(\beta,\frak a) = 0$ for all but finitely many $k$. The definition ensures that the form $h_\beta$ is nondegenerate on $\scr U^k(\beta,\frak a)$. It follows that $\scr J^k(\beta,\frak a)$ can be expressed as an {\it orthogonal\/} sum 
$$ 
\scr J^k(\beta,\frak a) = \sum_{i\ge k} \scr U^i(\beta,\frak a), 
$$ 
in which only finitely many terms are nonzero. 
\proclaim{Proposition} 
Let $\beta\in \Swr{\frak a}$ and set $E = F[\beta]$. There exists a unique $\frak o_F$-lattice $\frak I^1(\beta,\frak a)$ with the following properties:
\roster 
\item 
$\frak H^1(\beta,\frak a)\i \frak I^1(\beta,\frak a) \i \frak J^1(\beta,\frak a)$; 
\item 
$\frak I^1(\beta,\frak a)$ is stable under conjugation by $\r J(\beta,\frak a)$; 
\item 
the image $\scr I^1(\beta,\frak a)$ of\/ $\frak I^1(\beta,\frak a)$ in the alternating space $\scr J^1(\beta,\frak a)$ is a max\-imal totally isotropic subspace that is the sum of its intersections with the subspaces $\scr U^k(\beta,\frak a)$, $k\ge 1$. 
\endroster 
The lattice $\frak I^1(\beta,\frak a)$ has the following additional properties. 
\roster 
\item"\rm (4)" 
If $\beta'\in \Swr{\frak a}$ and $\scr C(\frak a,\beta',\psi^F) = \scr C(\frak a,\beta,\psi^F)$, then $\frak I^1(\beta',\frak a) =  \frak I^1(\beta,\frak a)$. 
\item"\rm (5)" 
If $K/F$ is a finite unramified extension, then 
$$ 
\frak I^1(\beta,\frak a_K) = \frak I^1(\beta,\frak a) \otimes_{\frak o_F}\frak o_K. 
\tag 1.4.2
$$ 
\endroster 
\endproclaim 
\demo{Proof} 
In the case of $G \cong \GL nF$, the result is 6.4 Proposition of \cite{5}. To deal with the first assertion in the general case, it is enough to show that there is a unique $\r J(\beta,\frak a)$-stable subspace $\scr I^1(\beta,\frak a)$ of $\scr J^1(\beta,\frak a)$ satisfying condition (3). 
\par 
Let $K/F$ be a finite unramified extension such that $A_K \cong \M nK$. The group $\scr J^1(\beta,\frak a)$ is a $\Bbbk_F$-vector space. Likewise, $\scr J^1(\beta,\frak a_K)$ is a $\Bbbk_K$-vector space and 1.3 Proposition implies 
$$ 
\scr J^k(\beta,\frak a_K) = \scr J^k(\beta,\frak a)\otimes_{\Bbbk_F}\Bbbk_K, \quad k\ge1. 
$$ 
Equally, if $\vG = \Gal KF$ then 
$$  
\scr J^k(\beta,\frak a) = \scr J^k(\beta,\frak a_K)^\vG, \quad k\ge1. 
$$ 
Indeed, $\scr V\mapsto \scr V^\vG$ is a bijection between the set of $\vG$-stable $\Bbbk_K$-subspaces $\scr V$ of $\scr J^1(\beta,\frak a_K)$ and the set of $\Bbbk_F$-subspaces of $\scr J^1(\beta,\frak a)$, the inverse being $\scr W\mapsto \scr W\otimes \Bbbk_K$. 
\par 
Remark also that, by the preceding lemma and the choice of $\psi^K$, the pairing $h^F_\beta$ is the restriction of $h^K_\beta$ to $\scr J^1(\beta,\frak a)$. 
\par 
We first prove that $\scr U^k(\beta,\frak a_K)$ is $\vG$-invariant. For $\gamma\in \vG$, let $\upr\gamma\psi^K$ be the character $x\mapsto \psi^K(x^\gamma)$, $x\in K$. There is a unique $t_\gamma\in U_K$ such that $\upr\gamma\psi^K(x) = \psi^K(t_\gamma x)$, $x\in K$. If $j\in \scr J^{1+k}(\beta,\frak a_K)$ and $u\in \scr U^k(\beta,\frak a_K)$, then 
$$ 
\align 
h_\beta^K(j,u^{\gamma}) &= \psi^K(\roman{tr}_{A_K}(\beta(u^{\gamma}j-ju^{\gamma}))) \\ 
&= \upr\gamma\psi^K(\roman{tr}_{A_K}(\beta(uj^{\gamma^{-1}}-j^{\gamma^{-1}}u)))  
= h_\beta^K(t_\gamma j^{\gamma^{-1}},u) = 1, 
\endalign 
$$ 
since $\scr J^{1+k}(\beta,\frak a_K)$ is $\vG$-invariant. Thus $u^\gamma \in \scr U^k(\beta,\frak a_K)$, as desired. It follows that 
$$ 
\scr J^k(\beta,\frak a) = \scr J^{1+k}(\beta,\frak a) \perp \scr U^k(\beta,\frak a_K)^\vG. 
$$ 
The summands here are $h_\beta^K$-orthogonal, so they are $h^F_\beta$-orthogonal whence 
$$ 
\scr U^k(\beta,\frak a_K)^\vG = \scr U^k(\beta,\frak a). 
$$ 
The uniqueness property of $\frak I^1(\beta,\frak a_K)$ implies that $\scr I^1(\beta,\frak a_K)$ is $\vG$-stable. Consider the subspace $\scr I^1(\beta,\frak a_K)^\vG$ of $\scr J^1(\beta,\frak a)$. It is the sum of its intersections with the spaces $\scr U^k(\beta,\frak a)$ and is totally isotropic. Comparing dimensions, it is a maximal totally isotropic subspace of $\scr J^1(\beta,\frak a)$. As $\scr I^1(\beta,\frak a_K)$ is stable under conjugation by $\r J(\beta,\frak a_K)$, so $\scr I^1(\beta,\frak a_K)^\vG$ is stable under conjugation by $\r J(\beta,\frak a_K)\cap G = \r J(\beta,\frak a)$. Thus $\scr I^1(\beta,\frak a_K)^\vG$ has all the properties demanded of $\scr I^1(\beta,\frak a)$. 
\par 
It remains to show that these properties determine $\scr I^1(\beta,\frak a) = \scr I^1(\beta,\frak a_K)^\vG$ uniquely. Let $\scr I_0$ be a maximal totally isotropic subspace of $\scr J^1(\beta,\frak a)$ satisfying the required conditions. The subspace $\scr I_0\otimes \Bbbk_K$ of $\scr J^1(\beta,\frak a_K)$ then has the necessary intersection property in (3). We show it is totally isotropic. Suppose the contrary. There then exist $x,y\in \scr I_0$ and a root of unity $\zeta$ in $K$ such that $h_\beta^K(1{+}x,1{+}\zeta y) \neq 1$. That is, 
$$ 
\align 
h_\beta^K(1{+}x,1{+}\zeta y) &= \psi^K(\roman{tr}_{A_K}(\beta\zeta(yx{-}xy))) \\ 
&= \psi^K(\zeta\roman{tr}_{A_K}(\beta(yx{-}xy))) \neq 1. 
\endalign 
$$ 
Therefore $\roman{tr}_{A_K}(\beta(yx{-}xy)) = \roman{tr}_A(\beta(yx{-}xy))$ does not lie in $\frak p_F$. Consequently, there exists $\zeta_0 \in \frak o_F$ such that 
$$ 
\psi^F(\zeta_0\roman{tr}_A(\beta(yx{-}xy))) = h^F_\beta(1{+}\zeta_0x,1{+}y) \neq 1. 
$$ 
That is, the space $\scr I_0$ is not totally isotropic. This contradiction implies that $\scr I_0\otimes \Bbbk_K$ is totally isotropic. It surely has properties (2) and (3), so it equals $\scr I^1(\beta,\frak a_K)$. Therefore $\scr I_0 = \scr I^1(\beta,\frak a_K)^\vG$, as required. 
\par 
If we fix $\beta$ for the moment, the construction of $\frak I^1(\beta,\frak a)$ has been done entirely in terms of a randomly chosen $\theta\in \scr C(\frak a, \beta,\psi^F)$. If $\theta$ also lies in $\scr C(\frak a,\beta',\psi^F)$, that is, if $\scr C(\frak a,\beta,\psi^F) = \scr C(\frak a,\beta',\psi^F)$, the uniqueness property of the first part implies $\frak I^1(\beta',\frak a) = \frak I^1(\beta,\frak a)$, as required for (4). Part (5) has already been done. \qed 
\enddemo 
Write 
$$ 
I^1(\beta,\frak a) = 1+\frak I^1(\beta,\frak a). 
\tag 1.4.3 
$$ 
\proclaim{Corollary} 
Let $f:A\to A'$ be an isomorphism of $F$-algebras. If $\frak a' = f(\frak a)$, then $\Swr{\frak a'} = f(\Swr{\frak a})$ and $I^1(f(\beta),\frak a') = f(I^1(\beta,\frak a))$, $\beta\in \Swr{\frak a}$. 
\endproclaim 
\demo{Proof} 
The first statement is 1.1 Proposition and the second follows from the uniqueness property of $\frak I^1(\beta,\frak a)$. \qed 
\enddemo 
The group $I^1(\beta,\frak a)$ is the {\it Lagrangian subgroup\/} of the section title and \cite{5}. 
\head 
2. Extensions of simple characters 
\endhead 
The notation is carried over from 1.4. We admit the case $p=2$ as far as the end of 2.1. 
\subhead 
2.1  
\endsubhead 
By definition, the group $I^1(\beta,\frak a)/H^1(\beta,\frak a)$ is a totally isotropic subspace of $J^1(\beta,\frak a)/H^1(\beta,\frak a)$, so a simple character $\theta\in \scr C(\frak a,\beta, \psi^F)$ admits extension to a character $\xi$ of the group $I^1(\beta,\frak a)$. The group $J^1(\beta,\frak a)/I^1(\beta,\frak a)$ acts on the set of such extensions by conjugation. Since $I^1(\beta,\frak a)/H^1(\beta,\frak a)$ is a maximal totally isotropic subspace of $J^1(\beta,\frak a)/H^1(\beta,\frak a)$, this set of extensions is a principal homogeneous space over  $J^1(\beta,\frak a)/I^1(\beta,\frak a)$. 
\proclaim{Lemma} 
Let $\xi$ be a character of $I^1(\beta,\frak a)$ such that $\theta = \xi\Mid H^1(\beta,\frak a)$ lies in $\scr C(\frak a,\beta,\psi^F)$. An element $g$ of $G$ intertwines $\xi$ if and only if $g\in \r J(\beta,\frak a)$ and $\xi^g = \xi$. 
\endproclaim 
\demo{Proof} 
If $g$ intertwines $\xi$, it surely intertwines $\theta$ and so lies in $\r J(\beta,\frak a)$. In particular, $g$ normalizes the character $\theta$ and also the group $I^1(\beta,\frak a)$. The lemma follows. \qed 
\enddemo 
In other words, the $G$-intertwining of $\xi$ is the $\r J(\beta,\frak a)$-normalizer of $\xi$. Our aim is to control this normalizer. To this end, we introduce a finer version of the set $\scr C(\frak a,\beta,\psi^F)$ of simple characters attached to $\beta\in \Swr{\frak a}$, following section 8.1 of \cite{5}. 
\definition{Definition} 
Let $\beta\in \Swr{\frak a}$ and define the positive integer $l$ by $\beta\frak a = \frak p^{-l}$. Let $\theta\in \scr C(\frak a,\beta,\psi^F)$. Say that $\theta$ is {\it adapted to $\beta$} if the following conditions hold. 
\roster 
\item If $l$ is even and a  jump of $\beta$, then $\theta\Mid H^{l/2}(\beta,\frak a) = \psi^A_\beta$. 
\item If $2k$ is a jump of $\beta$, such that $0<2k<l$, and if $[\frak a,l,2k,\gamma]$ is a simple stratum equivalent to $[\frak a,l,2k,\beta]$, there exists $\phi\in \scr C(\frak a,\gamma,\psi^F)$ such that 
$$ 
\theta\Mid H^k(\beta,\frak a) = \phi\,\psi^A_{\beta{-}\gamma}\Mid H^k(\beta,\frak a). 
$$ 
\endroster 
Let $\aC(\frak a,\beta,\psi^F)$ be the set of $\theta\in \scr C(\frak a,\beta,\psi^F)$ that are adapted to $\beta$. 
\enddefinition 
Note that condition (2) is independent of the choice of $\gamma$, {\it loc\. cit.} 
\proclaim{Proposition} 
Let $\beta\in \Swr{\frak a}$. 
\roster 
\item 
There exists $\theta\in \scr C(\frak a,\beta,\psi^F)$ that is adapted to $\beta$: the set $\aC(\frak a,\beta,\psi^F)$ is not empty. 
\item 
If $\vt\in \scr C(\frak a,\beta,\psi^F)$, there exists $\beta'\in \Swr{\frak a}$ such that $\vt\in \aC(\frak a,\beta',\psi^F)$. In particular, $ \scr C(\frak a,\beta',\psi^F) =  \scr C(\frak a,\beta,\psi^F)$ and $I^1(\beta',\frak a) =  I^1(\beta,\frak a)$. 
\endroster 
\endproclaim 
\demo{Proof} 
If $A$ is either a full matrix algebra or a division algebra, the result is proved in \cite{5} 8.1. The general case is identical and there is no need to repeat the details. \qed 
\enddemo 
\subhead 
2.2 
\endsubhead  
{\it From this point on, we assume that $p\neq 2$:\/} the immediate reason for this restriction is given in the Remark below. 
\definition{Definition} 
Let $\beta\in \Swr{\frak a}$, write $E = F[\beta]$ and set 
$$ 
\r I(\beta,\frak a) = E^\times I^1(\beta,\frak a). 
\tag 2.2.1 
$$ 
Define $\scr D(\frak a,\beta,\psi^F)$ to be the set of characters $\lambda$ of the group $\r I(\beta,\frak a)$ such that $\lambda\Mid H^1(\beta,\frak a) \in \aC(\frak a,\beta,\psi^F)$. Set 
$$ 
\scr D(\frak a,\psi^F) = \bigcup_{\beta\in \Swr{\frak a}} \scr D(\frak a,\beta,\psi^F). 
\tag 2.2.2 
$$ 
\enddefinition 
\proclaim{Proposition} 
Let $\beta\in \Swr{\frak a}$ and write $E = F[\beta]$. 
\roster 
\item 
The restriction map $\scr D(\frak a,\beta,\psi^F) \to \aC(\frak a,\beta,\psi^F)$ is surjective. 
\item 
Let $K/F$ be a finite unramified extension. If $\lambda\in \scr D(\frak a_K,\beta,\psi^K)$, the character 
$$
\lambda^F = \lambda\Mid \r I(\beta,\frak a) 
\tag 2.2.3 
$$ 
lies in $\scr D(\frak a,\beta,\psi^F)$. The map 
$$ 
\align 
\scr D(\frak a_K,\beta,\psi^K) &\longrightarrow \scr D(\frak a,\beta,\psi^F), \\ 
\lambda &\longmapsto \lambda^F, 
\endalign 
$$ 
is surjective. 
\item 
If $f:A\to A'$ is an isomorphism of $F$-algebras, then $f(\r I(\beta,\frak a)) = \r I(f(\beta), f(\frak a))$ and $f$ induces a bijection 
$$ 
\align 
\scr D(\frak a,\beta,\psi^F) &\longrightarrow \scr D(f(\frak a), f(\beta),\psi^F), \\ 
\lambda &\longmapsto \lambda\circ f^{-1}. 
\endalign 
$$ 
\endroster 
\endproclaim 
\demo{Proof} 
Suppose first that $A\cong \M nF$. Let $\theta \in \aC(\frak a,\beta,\psi^F)$. As in 8.4 Proposition of \cite{5}, there is a character $\xi$ of $I^1(\beta,\frak a)$ that extends $\theta$ and is stable under conjugation by $E^\times$. The character $\xi$ then extends to a character $\lambda$ of $\r I(\beta,\frak a)$, whence follows (1) in this case. 
\par 
In the general case, let $K/F$ be unramified of degree divisible by $n$, so that $A_K \cong \M nK$. 
\proclaim{Lemma 1} 
The restriction map $\scr C(\frak a_K,\beta,\psi^K) \to \scr C(\frak a,\beta,\psi^F)$ induces a surjection   $\aC(\frak a_K,\beta,\psi^K) \to \aC(\frak a,\beta,\psi^F)$. 
\endproclaim 
\demo{Proof} 
If $\vt \in \aC(\frak a_K,\beta,\psi^K)$, then 2.1 Definition implies that $\vt\Mid H^1(\beta,\frak a)$ lies in $\aC(\frak a,\beta,\psi^F)$. So, we take $\theta\in \aC(\frak a,\beta,\psi^F)$ and construct a character $\theta_K \in \aC(\frak a_K,\beta,\psi^K)$ that extends $\theta$. 
\par 
Let $t\ge0$ be the least integer for which there exists $\vt\in \scr C(\frak a_K,\beta,\psi^K)$ agreeing with $\theta$ on $H^{1+t}(\beta,\frak a)$ and satisfying the conditions of the definition relative to any even jump $2k$, with $k > t$. 
\par 
If $l$ is even and a jump of $\beta$, the definition yields $t<l/2$. If $l$ is the only even jump, we get $t=0$ since there are no further restrictions to be observed. In all other cases, we still have $t<l/2$. If $t<k$ for any even jump $2k$, we again get $t=0$. Otherwise, let $2k$ be the greatest even jump such that $k\le t$. We may adjust our choice of $\vt$ so that it agrees with $\theta$ on $H^{1+k}(\beta,\frak a)$ without affecting the conditions already imposed in our hypothesis. That is, we can assume that $t = k$ where $2k$ is an even jump. Take a simple stratum $[\frak a,l,2k,\gamma]$ equivalent to $[\frak a,l,m,\beta]$. On $H^k(\beta,\frak a)$, $\theta$ takes the form $\phi\psi^A_{\beta{-}\gamma}$, for some $\phi\in \scr C(\frak a,\gamma,\psi^F)$. The standard construction of simple characters implies that, on $H^{1+k}(\beta,\frak a_K)$, the character $\vt$ takes the form $\vf\psi^{A_K}_{\beta{-}\gamma}$, for some $\vf\in \scr C(\frak a_K,\gamma,\psi^K)$. Surely $\vf$ agrees with $\phi$ on $H^{1+k}(\beta,\frak a)$. However, 7.1 Theorem of \cite{5} implies that $\vf$ agrees with $\phi$ on $H^k(\beta,\frak a)$. We could therefore have chosen our original $\vt$ to agree with $\theta$ on the larger group $H^k(\beta,\frak a)$, contrary to our definition of $t$. So, in all cases, $t=0$ and the lemma is proven. \qed
\enddemo 
\remark{Note} 
Lemma 1 and its proof remain valid when $p=2$. 
\endremark 
Continuing with the proof of the proposition, let $\theta\in \aC(\frak a,\beta,\psi^F)$. Take $\theta_K \in \aC(\frak a_K,\beta,\psi^K)$ such that $\theta = \theta_K\Mid H^1(\beta,\frak a)$. Since $A \cong \M nK$, we know that $\theta_K$ admits extension to a character $\lambda_K$ of $\r I(\beta,\frak a_K)$. The restriction $\lambda^F_K = \lambda_K\Mid \r I(\beta,\frak a)$ provides an extension of $\theta$, and part (1) of the proposition is proven. 
\par 
We show that any extension of $\theta$ to $\r I(\beta,\frak a)$ arises in this way. As a first step, let $\vt$ be a character of $I^1(\beta,\frak a)$ that extends $\theta$ and is stable under conjugation by $E^\times$. It differs from $\lambda_K\Mid I^1(\beta,\frak a)$ by a character $\phi$ of $I^1(\beta,\frak a)/H^1(\beta,\frak a)$ stable under conjugation by a prime element $\vp$ of $E$. 
\proclaim{Lemma 2} 
The character $\phi$ extends to a character $\phi_K$ of $I^1(\beta,\frak a_K)/H^1(\beta,\frak a_K)$ stable under conjugation by $\vp$. 
\endproclaim 
\demo{Proof} 
Identify $I^1(\beta,\frak a)/H^1(\beta,\frak a)$ with the $\Bbbk_F$-space $\frak I^1(\beta,\frak a)/\frak H^1(\beta,\frak a)$. Viewed as a character of this group, $\phi$ is trivial on the image of the map $A_\vp:x\mapsto \vp x\vp^{-1}{-}x$. As $I^1(\beta,\frak a_K)/H^1(\beta,\frak a_K) = I^1(\beta,\frak a)/H^1(\beta,\frak a)\otimes_{\Bbbk_F} \Bbbk_K$, the result follows straightaway. \qed 
\enddemo 
The character $\phi_K\otimes \lambda_K\Mid I^1(\beta,\frak a_K)$ is stable under conjugation by $\vp$, and hence by $KE^\times$. It therefore extends to a character of $\r I(\beta,\frak a_K)$ and some such extension restricts to $\lambda$. This proves (2) in the case where $[K{:}F]$ is divisible by $n$, and the general case follows by transitivity. 
\par 
In part (3), the map $f$ carries $\aC(\frak a,\beta,\psi^F)$ bijectively to $\aC(f(\frak a), f(\beta),\psi^F)$, as follows directly from the definition. \qed 
\enddemo 
\remark{Remark} 
The definition here of the set $\scr D(\frak a,\beta,\psi^F)$ is different from, and more inclusive than, the one used in \cite{5}. We have found it more convenient. When $p=2$, both approaches fail: there are examples where $\scr D(\frak a,\beta,\psi^F)$ is empty (\cite{5} 8.3 Remark). 
\endremark 
\subhead 
2.3 
\endsubhead 
We relate the sets $\scr D(\frak a,\beta,\psi^F)$, $\scr D(\frak a,\psi^F)$  to representations of $G$. Observe that, as a consequence of part (3) of 2.2 Proposition, the group $\r K_\frak a$ acts on the set $\scr D(\frak a,\psi^F)$ by conjugation. 
\proclaim{Theorem} 
Let $\frak a$ be a minimal hereditary $\frak o_F$-order in $A$, and let $\beta\in \Swr{\frak a}$. 
\roster 
\item 
If $\lambda\in \scr D(\frak a,\beta,\psi^F)$, the induced representation 
$$ 
\pi_G(\lambda) = \cind_{\,\r I(\beta,\frak a)}^G\,\lambda 
$$ 
is irreducible and cuspidal. Its equivalence class lies in $\awr G$. 
\item  
The map 
$$ 
\lambda \longmapsto \pi_G(\lambda), \quad \lambda\in \scr D(\frak a,\psi^F), 
$$ 
induces a canonical bijection 
$$ 
\pi_G: \r K_\frak a\backslash \scr D(\frak a,\psi^F) @>{\ \ \approx\ \ }>> \awr G. 
\tag 2.3.1 
$$ 
\endroster 
\endproclaim 
\demo{Proof} 
Let $\lambda \in \scr D(\frak a,\beta,\psi^F)$, and put $\vt = \lambda\Mid I^1(\beta,\frak a)$, $\theta = \lambda\Mid H^1(\beta,\frak a)$. There is a unique irreducible representation $\eta$ of $J^1(\beta,\frak a)$ containing $\theta$ \cite{18} 2.2. By definition, $I^1(\beta,\frak a)/H^1(\beta,\frak a)$ is a maximal isotropic subspace of the alternating space $J^1(\beta,\frak a)/H^1(\beta,\frak a)$, so $\eta$ is induced by any character of $I^1(\beta,\frak a)$ extending $\theta$. If $\vL$ is the representation of $\r J(\beta,\frak a) $ induced by $\lambda$, the Mackey restriction formula shows that $\vL\Mid J^1(\beta,\frak a)$ is the irreducible representation $\eta$. Therefore $\vL$ is irreducible. Any $g\in G$ that intertwines $\theta$ lies in $\r J(\beta,\frak a)$, so 
$$ 
\pi_G(\lambda) = \cind_{\r I(\beta,\frak a)}^G\,\lambda = \cind_{\r J(\beta,\frak a)}^G\,\vL 
\tag 2.3.2 
$$ 
is irreducible and cuspidal. Since it contains the simple character $\theta\in \scr C(\frak a,\beta,\psi^F)$ and the field extension $F[\beta]/F$ is totally ramified of degree $n$, the representation $\pi_G(\lambda)$ is totally ramified of parametric degree $n$. That is, $\pi_G(\lambda) \in \awr G$. 
\par 
Conversely, let $\pi \in \awr G$. As in \cite{20}, there exists an extended maximal simple type $\vL$ in $G$, inducing $\pi$. Since $\pi \in \awr G$, $\vL$ is a representation of a group $\r J(\beta,\frak a)$, for some $\beta\in \Swr{\frak a}$, that contains a simple character $\theta\in \scr C(\frak a,\beta,\psi^F)$. By 2.1 Proposition, we may assume $\theta\in \aC(\frak a,\beta,\psi^F)$. Let $\vt$ be a character of $I^1(\beta,\frak a)$, extending $\theta$. Any two choices of $\vt$ are $J^1(\beta,\frak a)$-conjugate, so $\vt$ occurs in $\vL$. By part (1) of 2.2 Proposition, we may take $\vt$ to be $F[\beta]^\times$-stable, so there exists $\lambda\in \scr D(\frak a,\beta,\psi^F)$ extending $\vt$ and occurring in $\vL$. The representation of $\r J(\beta,\frak a)$ induced by $\lambda$ is then $\vL$, giving $\pi = \pi_G(\lambda)$, as desired. 
\par 
The group $\r K_\frak a$  acts on the set $\scr D(\frak a,\psi^F)$ by conjugation. For $\lambda\in \scr D(\frak a,\psi^F)$, the equivalence class of $\pi_G(\lambda)$ depends only on the $\r K_\frak a$-orbit of $\lambda$, so $\lambda \mapsto \pi_G(\lambda)$ induces a surjective map $\r K_\frak a\backslash \scr D(\frak a,\psi^F) \to \awr G$. To prove it is injective, take $\beta_i \in \Swr{\frak a}$ and $\lambda_i\in \scr D(\frak a,\beta_i,\psi^F)$, $i=1,2$, and suppose that $\pi_G(\lambda_1) = \pi_G(\lambda_2) =\pi$, say. The simple characters $\theta_i = \lambda_i\Mid H^1(\beta_i,\frak a)$ intertwine in $G$. They are therefore $\r K_\frak a$-conjugate \cite{20} Theorem 6.1. We may assume they are equal, say $\theta_1 = \theta_2 = \theta$, implying that the sets $\scr C(\frak a,\beta_i,\psi^F)$ are the same. In particular, $I^1(\beta_1,\frak a) = I^1(\beta_2,\frak a)$ (1.4 Proposition) and the same holds for the $H^1$, $J^1$ and $\r J$ groups. After applying a $J^1(\beta_i,\frak a)$-conjugation, we can assume that the $\lambda_i$ agree on $I^1(\beta_i,\frak a)$ and intertwine in $G$. This intertwining is implemented by an element $x$ which intertwines $\theta$ and so lies in $\r J(\beta_i,\frak a)$. The element $x$ normalizes $I^1(\beta_i,\frak a)$ and fixes the character $\lambda_i\Mid I^1(\beta_i,\frak a)$. Therefore $x\in \r I(\beta_i,\frak a)$ and $x$ fixes $\lambda_i$. Therefore $\lambda_1 = \lambda_2$, as required. \qed 
\enddemo 
\subhead 
2.4 
\endsubhead  
Let $\beta\in \Swr{\frak a}$, let $\theta\in \aC(\frak a,\beta,\psi^F)$ and write $E = F[\beta]$. Following \cite{5} 8.3, 8.4, we write down an $E^\times$-invariant character $\pre I\theta$ of $I^1(\beta,\frak a)$. 
\par 
Let $\beta^{-1}\frak a = \frak p^l$, where $\frak p$ is the Jacobson radical of $\frak a$. If $\beta$ has no even jumps, that is, if $H^1(\beta,\frak a) = J^1(\beta,\frak a)$, then $I^1(\beta,\frak a) = H^1(\beta,\frak a)$ and we set $\pre I\theta = \theta$. This certainly satisfies the requirement. Otherwise, let $2s_1<2s_2< \dots < 2s_t$ be the even jumps of $\beta$. Reverting to the notation of 1.4, let $\scr I$ be the image of $I^1(\beta,\frak a)$ in $J^1(\beta,\frak a)/H^1(\beta,\frak a)$ and let $I_k$ be the inverse image, in $J^{s_k}(\beta,\frak a)$ of $\scr I\cap \scr U^k(\beta,\frak a)$, $1\le k\le t$. We then have 
$$ 
I^1(\beta,\frak a) = H^1(\beta,\frak a)\,I_1\,I_2\, \dots\,I_t, 
\tag 2.4.1 
$$ 
with all factors in the product commuting modulo the kernel of $\theta$. If $s_{2t} = l$, we define 
$$ 
\pre I\theta(1{+}x) = \psi^A_\beta(1{+}x-x^2/2), \quad 1{+}x \in I_t. 
\tag 2.4.2 
$$ 
Otherwise, let $s_{2k} < l$ and choose a simple stratum $[\frak a,l,s_{2k},\gamma_j]$ equivalent to $[\frak a,l,s_{2k},\beta]$. By 2.1 Definition, $\theta\Mid H^{s_k}(\beta,\frak a) = \vt\psi^A_{\beta{-}\gamma_k}\Mid H^{s_k}(\beta,\frak a)$, for some $\vt\in \scr C(\frak a,\gamma_k,\psi^F)$. We set 
$$
\pre I\theta(1{+}x) = \vt(1{+}x)\,\psi^A_{\beta{-}\gamma_k}(1{+}x{-}x^2/2), \quad 1{+}x\in I_k. 
\tag 2.4.3 
$$ 
Taking into account the product formula (2.4.1), the expressions (2.4.2), (2.4.3) define $\pre I\theta$ as a character of $I^1(\beta,\frak a)$. As in \cite{5}, the character $\pre I\theta$ is stable under conjugation by $E^\times$. 
\par 
To go a step further, note that $E^\times \cap H^1(\beta,\frak a) =  E^\times \cap I^1(\beta,\frak a) = U^1_E$. Let $\xi$ be a character of $E^\times$ agreeing with $\theta$ on $U^1_E = E^\times \cap H^1(\beta,\frak a)$. The formula 
$$ 
\xi\odot\theta:ux \longmapsto \xi(u)\,\pre I\theta(x),\quad u\in U^1_E,\ x\in I^1(\beta,\frak a), 
\tag 2.4.4 
$$ 
defines $\xi\odot \theta$ as a character of $\r I(\beta,\frak a)$. Surely, $\xi\odot\theta \in \scr D(\frak a,\beta,\psi^F)$. 
\proclaim{Proposition} 
Let $\pi\in \awr G$ contain the simple character $\theta\in\aC(\frak a,\beta,\psi^F)$. Set $E = F[\beta]$. 
\roster 
\item 
The representation $\pi$ contains the character $\pre I\theta$ of $I^1(\beta,\frak a)$. 
\item 
There is a unique character $\xi$ of $E^\times$, agreeing with $\theta$ on $U^1_E$, such that $\pi$ contains $\xi\odot \theta$ and, consequently, $\pi = \pi_G(\xi\odot\theta)$. 
\endroster 
The map $\xi\mapsto \pi_G(\xi\odot\theta)$ is a bijection between the set of characters $\xi$ of $E^\times$, that agree with $\theta$ on $U^1_E$, and the set of elements of $\awr G$ that contain $\theta$. 
\endproclaim 
\demo{Proof} 
Surely $\pi$ contains a character $\phi$ of $I^1(\beta,\frak a)$ extending $\theta$. Since the space $I^1(\beta,\frak a)/H^1(\beta,\frak a)$ is a maximal totally isotropic subspace of $J^1(\beta,\frak a)/H^1(\beta,\frak a)$, the character $\phi$ is $J^1(\beta,\frak a)$-conjugate to $\pre I\theta$, which therefore occurs in $\pi$. If $\eta$ is the unique irreducible representation of $J^1(\beta,\frak a)$ that contains $\theta$, then $\eta$ occurs in $\pi$ with multiplicity one. So, there is a unique character $\lambda\in \scr D(\frak a,\beta,\psi^F)$ that occurs in $\pi$ and extends $\pre I\theta$. Surely there exists a unique character $\xi$, of the required form, such that $\lambda = \xi\odot \theta$. \qed 
\enddemo 
\subhead 
2.5 
\endsubhead 
This construction behaves properly with respect to unramified base field extension. 
\proclaim{Proposition} 
Let $\beta\in \Swr{\frak a}$, let $\theta\in \aC(\frak a,\beta,\psi^F)$ and let $\xi$ be a character of $E^\times = F[\beta]^\times$ agreeing with $\theta$ on $U^1_E$. Let $K/F$ be a finite unramified extension. If $\theta_K\in \aC(\frak a_K,\beta,\psi^K)$ agrees with $\theta$ on $H^1(\beta,\frak a)$, and if $\xi_K$ is a character of $KE^\times$ agreeing with $\theta_K$ on $U^1_{KE}$ and with $\xi$ on $E^\times$, then 
 $$ 
\align 
\pre I\theta_K\Mid I^1(\beta,\frak a) &= \pre I\theta, \\
\xi_K\odot \pre I\theta_K\Mid \r I(\beta,\frak a) &= \xi\odot\pre I\theta. 
\endalign  
$$ 
\endproclaim 
\demo{Proof} 
Immediate. \qed 
\enddemo 
\head 
3. Change of group and endo-classes
\endhead 
Let $A$ and $B$ be central simple $F$-algebras of dimension $n^2$, where $n$ is a power of $p$. Let $\frak a$ and $\frak b$ be minimal hereditary $\frak o_F$-orders in $A$ and $B$ respectively. Write $G = A^\times$ and $H = B^\times$. Other notation is carried over from sections 1 and 2. 
\par 
If $K/F$ is an unramified extension of finite degree divisible by $n$, the $K$-algebras $A_K$ and $B_K$ are both isomorphic to $\M nK$. In this section, we use such isomorphisms to exploit the naturality properties of the sets $\scr D(\frak a,\beta,\psi^F)$ laid out in 2.2 Proposition. In the case to hand, these properties restrict to a transfer of simple characters between $G$ and $H$ that preserves endo-classes. 
\par 
More to the point, we obtain a process, called {\it parametric transfer,\/} for moving representations between $\awr G$ and $\awr H$. By the end of the section, it still depends on one choice made in the construction.  
\subhead 
3.1 
\endsubhead 
We start with a basic formal result. 
\proclaim{Proposition} 
Let $\beta\in \Swr{\frak a}$ and write $E = F[\beta]$. Let $K/F$ be a finite unramified extension of degree divisible by $n$. 
\roster 
\item 
There exists an $F$-embedding $\iota:E\to B$ such that $\iota(E^\times) \i \r K_\frak b$ and $\iota(\beta) \in \Swr{\frak b}$. An $F$-embedding $\iota':E\to B$ has the same property if and only if $\iota' = \roman{Ad}\,x\circ \iota$, for some $x\in \r K_\frak b$. 
\item 
The map $\iota$ extends to an isomorphism $\iota_K:A_K \to B_K$ of $K$-algebras such that $\iota_K(\frak a_K) = \frak b_K$. Any such extension $\iota_K$ has the property
$$ 
\aligned 
\iota_K(I^1(\beta,\frak a_K)) &= I^1(\iota(\beta),\frak b_K), \\ 
\iota_K(\r I(\beta,\frak a_K)) & = \r I(\iota(\beta), \frak b_K), 
\endaligned 
\tag 3.1.1 
$$ 
and induces a bijection 
$$ 
\aligned 
\scr D(\frak a_K,\beta,\psi^K) &\longrightarrow \scr D(\frak b_K,\iota(\beta),\psi^K), \\ 
\lambda &\longmapsto \lambda\circ \iota_K^{-1}, 
\endaligned 
\tag 3.1.2
$$ 
\item 
If $\iota'_K:A_K \to B_K$ also extends $\iota$ and has the property $\iota'_K(\frak a_K) = \frak b_K$, there exists $y\in KE^\times$ such that $\iota'_K = \iota_K \circ \roman{Ad}\,y$. 
\endroster 
\endproclaim  
\demo{Proof} 
There surely exists an $F$-embedding $\iota: E\to B$. Since $E/F$ is totally ramified of degree $n$, there is a {\it unique\/} minimal hereditary $\frak o_F$-order $\frak b_1$ in $B$ such that $\iota(E^\times) \i \r K_{\frak b_1}$, as in 1.1. Replacing $\iota$ by $\roman{Ad}\,x\circ\iota$, for some $x\in H$, we may take $\frak b_1 = \frak b$. That is, $\iota$ satisfies the first assertion of (1) and the second follows from 1.1 Proposition. The third assertion of (1) follows from the uniqueness of $\frak b$. 
\par 
We extend $\iota$, by linearity, to a $K$-embedding $\iota_K$ of the field $K\otimes_FE = KE = K[\beta]$ in $B_K$. Since $\iota(\beta) \in \Swr{\frak b}$, we have $\iota_K(\beta) \in \Swr{\frak b_K}$ ({\it cf\.} 1.3). Since $A_K \cong B_K$, the Skolem--Noether theorem implies that $\iota_K$ extends to an isomorphism $A_K\to B_K$ of $K$-algebras. The image $\iota_K(\frak a_K)$ is a hereditary $\frak o_K$-order stable under conjugation by $\iota_K(K[\beta])^\times = K[\iota(\beta)]^\times$, so $\iota_K(\frak a_K) = \frak b_K$. The map $\iota_K$ carries $H^1(\beta,\frak a_K)$ to $H^1(\iota(\beta),\frak b_K)$ and similarly for the $J^1$-groups. Further, the map $\theta \mapsto \theta\circ \iota_K^{-1}$ is a bijection 
$$ 
\scr C(\frak a_K,\beta,\psi^K) \longrightarrow \scr C(\frak b_K,\iota(\beta),\psi^K). 
\tag 3.1.3 
$$ 
The property (3.1.1) follows from 2.2 Proposition (3). Moreover, the map (3.1.3) restricts to a bijection 
$$ 
\aC(\frak a_K,\beta,\psi^K) \longrightarrow \aC(\frak b_K,\iota(\beta),\psi^K), 
$$ 
by 2.2 Lemma 1. Thus (3.1.2) is a bijection as required. 
\par 
Finally, if $\iota'_K$ is another extension as in (3), then $\iota'_K = \iota_K\circ \roman{Ad}\,z$, for some $z\in G$. It also agrees with $\iota_K$ on $KE$, while the field $KE$ is its own centralizer in $A_K$.  \qed 
\enddemo 
Take an $F$-embedding $\iota:E\to B$ and an extension $\iota_K:A_K \to B_K$, as in the proposition. Composing the bijection $ \scr D(\frak a_K,\beta,\psi^K) \to \scr D(\frak b_K,\iota(\beta),\psi^K)$ of (3.1.2) with the restriction map $\scr D(\frak b_K,\iota(\beta),\psi^K) \to \scr D(\frak b,\iota(\beta),\psi^F)$ of 2.2 Proposition, we get a surjective map 
$$ 
\aligned 
\scr D(\frak a_K,\beta,\psi^K) &\longrightarrow \scr D(\frak b,\iota(\beta),\psi^F), \\
(\r I,\kappa) &\longmapsto (\r I\cap H,\iota\kappa^F) , 
\endaligned 
\tag 3.1.4 
$$ 
where we abbreviate $\r I = \r I(\iota(\beta),\frak b_K)$, so that $\r I\cap H = \r I(\iota(\beta, \frak b)$, and 
$$ 
\iota\kappa^F = \kappa\circ\iota^{-1}\Mid \r I\cap H. 
$$ 
We form the representation $\pi_H(\iota\kappa^F) =\cind_{\r I\cap H}^H\,\iota\kappa^F \in \awr H$, as in 2.3. 
\proclaim{Corollary} 
Let $\iota$, $\iota':E\to B$ be embeddings as in the proposition. If $\kappa \in \scr D(\frak a_K,\beta, \psi^K)$, then $\pi_H(\iota\kappa^F) \cong \pi_H(\iota'\kappa^F)$. 
\endproclaim 
\demo{Proof} 
Part (3) of the proposition shows that $\pi_H(\iota\kappa^F)$ depends only on $\iota$, not on the choice of extension $\iota_K$, while (1) shows it is independent of $\iota$. \qed 
\enddemo 
\subhead 
3.2 
\endsubhead 
The procedure of 3.1 is essentially independent of the choice of $K/F$. For, if $L/K$ is a finite unramified extension, we have a commutative diagram 
$$ 
\CD 
\scr D(\frak a_L,\beta,\psi^L) @>{\ \ \approx\ \ }>> \scr D(\frak b_L,\iota(\beta),\psi^L) \\ 
@VVV @VVV \\  
\scr D(\frak a_K,\beta,\psi^K) @>{\ \ \approx\ \ }>> \scr D(\frak b_K,\iota(\beta),\psi^K) 
\endCD \tag 3.2.1 
$$ 
in which the vertical arrows are the surjective restriction maps. 
\subhead 
3.3 
\endsubhead 
Take $\theta\in \aC(\frak a,\beta,\psi^F)$. It is the restriction of some $\theta_K \in \aC(\frak a_K,\beta,\psi^K)$ (2.2 Lemma 1). We use an embedding $\iota$, satisfying the conditions of 3.1 Proposition, to form 
$$ 
\theta_K\circ\iota_K^{-1} \in \aC(\frak b_K,\iota\beta,\psi^K). 
$$ 
We set 
$$ 
\iota\theta_K^F = \theta_K\circ\iota_K^{-1}\Mid H^1(\iota\beta,\frak b). 
\tag 3.3.1 
$$ 
Thus $\iota\theta_K^F \in \aC(\frak b,\iota\beta,\psi^F)$. Using the language of \cite{3}, we have: 
\proclaim{Proposition} 
The simple characters $\theta$, $\iota\theta_K^F$ are endo-equivalent. 
\endproclaim 
\demo{Proof} 
This is Theorem 1.13 plus Remark 6.9 of \cite{3}. \qed 
\enddemo 
\subhead 
3.4 
\endsubhead 
Let $\pi\in \awr G$. As in 2.3, there is a character $\lambda\in \scr D(\frak a,\psi^F)$ that induces $\pi$, and this $\lambda$ is unique up to $\r K_\frak a$-conjugation. We may choose $\beta\in \Swr{\frak a}$ so that $\lambda\in \scr D(\frak a,\beta,\psi^F)$. Following the procedure of 3.1, we set $E = F[\beta]$ and choose 
\roster 
\item 
an $F$-embedding $\iota:E\to B$ such that $\iota(E)^\times \i \r K_\frak b$, 
\item 
a finite unramified field extension $K/F$ of degree divisible by $n$, 
\item 
a $K$-isomorphism $\iota_K:A_K \to B_K$ extending $\iota$, and 
\item a character $\lambda_K\in \scr D(\frak a_K,\beta,\psi^K)$ extending $\lambda$. 
\endroster 
Having made these choices, we get a representation $\pi' = \pi_H(\iota\lambda^F_K) \in \awr H$. Following (3.1) Corollary, $\pi'$ actually depends only on the choice of $\lambda^K$ in (4). We make no effort at this stage to eliminate that dependence. We say that a representation $\pi'\in \awr H$, obtained from $\pi\in \awr G$ by such a choice, is a {\it parametric transfer of\/} $\pi$. 
\par 
Observe that, if we have a third algebra $C$, a representation $\pi''\in \awr {C^\times}$ that is a parametric transfer of $\pi'$ (relative to the same $\beta$) is also a parametric transfer of $\pi$. 
\subhead 
3.5 
\endsubhead 
We look back at the constructions of 2.4, 2.5. Thus we take $\beta\in \Swr{\frak a}$, $\theta\in \aC(\frak a,\beta,\psi^F)$ and set $E = F[\beta]$. Let $\xi$ be a character of $E^\times$ agreeing with $\theta$ on $U^1_E$. We form the character $\xi\odot\theta\in \scr D(\frak a,\beta,\psi^F)$ as in (2.4.4). Let $\theta_K \in \aC(\frak a_K,\beta,\psi^K)$ agree with $\theta$ on $H^1(\beta,\frak a)$. By 2.5 Proposition, $\pre I\theta = \pre I\theta_K\Mid H^1(\beta,\frak a)$. So, if $\xi_K$ is a character of $KE^\times$ that agrees with $\theta_K$ on $U^1_{KE}$ and with $\xi$ on $E^\times$, we have $\xi_K\odot \theta_K \in \scr D(\frak a_K,\beta,\psi^K)$ and 
$$ 
\xi_K\odot \theta_K \Mid \r I(\beta,\frak a) = \xi\odot\theta. 
\tag 3.5.1 
$$ 
Looking back to 3.4 and chasing through the definitions, we find: 
\proclaim{Proposition} 
If $\lambda = \xi\odot \theta$ and $\lambda_K = \xi_K\odot\theta_K$, then 
$$ 
\iota\lambda^F_K = (\xi\circ\iota^{-1})\odot \iota\theta_K^F. 
$$ 
\endproclaim 
\head 
4. Transfer via completion 
\endhead 
We analyze more deeply the embeddings $\iota$ of 3.1 using a technique of passing to a limit, as suggested by the diagram (3.2.1). 
\par 
Apart from results in 4.5 concerning the character sets $\scr D$, everything in this section holds equally when $p=2$. The notation follows on from the preceding sections but, from 4.3 onwards, it is convenient to choose our ``base point'' $A$ to be the matrix algebra $\M nF$. That entails no loss of generality. 
\subhead 
4.1 
\endsubhead 
We need some new notation. 
\par 
Let $F_\infty/F$ be a maximal unramified extension of $F$. Thus $F_\infty$ is the union of all finite unramified extensions  $K/F$ inside some algebraic closure of $F$. The discrete valuation $\ups_F$ on $F$ extends to a discrete valuation $F^\times_\infty \to \Bbb Z$, also denoted $\ups_F$. The associated discrete valuation ring is $\frak o_\infty = \bigcup \frak o_K$, with $K$ ranging as before. The maximal ideal of $\frak o_\infty$ is $\vp_F\frak o_\infty$, where $\vp_F$ is a prime element of $F$. The residue field $\Bbbk_\infty = \frak o_\infty/\frak p_\infty$ is an algebraic closure of the residue field $\Bbbk_F$ of $F$. 
\par 
Let $\wt F$ be the completion of $F_\infty$ with respect to $\ups_F$. Thus $\ups_F$ extends to a discrete valuation on $\wt F$. The associated discrete valuation ring is the closure of $\frak o_\infty$ in $\wt F$: we denote it by $\tilde{\frak o}$. This has maximal ideal $\vp_F\tilde{\frak o}$. The residue field $\wt\Bbbk = \tilde{\frak o}/\vp_F\tilde{\frak o}$ is equal to $\Bbbk_\infty$. The group $\wt F^\times$ is generated by $\vp_F$, the group $\tilde{\bs\mu}$ of roots of unity in $F_\infty$ of order relatively prime to $p$, and the principal unit group $1{+}\vp_F\tilde{\frak o}$. 
\par 
Let $\vO = \Gal{F_\infty}F$. Thus $\vO$ is procyclic and canonically isomorphic to $\Gal{\wt\Bbbk}{\Bbbk_F}$. It is topologically generated by the arithmetic Frobenius $\sigma_F$, that acts on $\tilde{\bs\mu}$ as $\zeta\mapsto \zeta^q$, where $q = |\Bbbk_F|$. Every element of $\vO$ extends uniquely to a continuous $F$-automorphism of $\wt F$, and $\vO$ is so identified with the group of continuous $F$-automorphisms of $\wt F$. If $K/F$ is a finite unramified extension and $\vO_K = \Gal{F_\infty}K$, the set of $\vO_K$-fixed points in $\wt F$ is again $K$. 
\subhead 
4.2 
\endsubhead 
We return to the situation of section 2. Thus $n = p^r$, for an integer $r\ge1$, and $A$ is a central simple $F$-algebra of dimension $n^2$. Let $\frak a$ be a minimal hereditary $\frak o_F$-order in $A$. Let $\beta\in \Swr{\frak a}$ and set 
$$ 
\r I_\infty(\beta,\frak a) = \bigcup_{K/F} \r I(\beta,\frak a_K). 
$$ 
Here $K/F$ ranges over the finite sub-extensions of $F_\infty/F$ and the union is taken in $\bigcup_K A_K = A\otimes F_\infty$. Let $\tilde{\r I}(\beta,\frak a)$ be the closure of $\r I_\infty(\beta,\frak a)$ in the topological group $G_{\wt F} = \big(A\otimes_F \wt F \big)^\times$. 
\proclaim{Proposition} 
Let $\beta\in \Swr{\frak a}$ and define 
$$ 
\widetilde{\frak I}^1(\beta,\frak a) = \frak I^1(\beta,\frak a)\otimes_{\frak o_F} \tilde{\frak o},\quad \tilde I^1(\beta,\frak a) = 1+\widetilde{\frak I}^1(\beta,\frak a). 
$$ 
We then have $\tilde{\r I}(\beta,\frak a) = \widetilde F[\beta]^\times \tilde I^1(\beta,\frak a)$. If $K/F$ is a finite sub-extension of $\widetilde F/F$, then 
$$ 
\align 
I^1(\beta,\frak a_K) &= \tilde I^1(\beta,\frak a)\cap G_K, \\ 
\r I(\beta,\frak a_K) &= \tilde{\r I}(\beta,\frak a)\cap G_K. 
\endalign 
$$ 
\endproclaim 
The proof is immediate. 
\definition{Definition}  
Choose, once for all, a character $\tilde\psi$ of $\widetilde F$, of level one, such that $\tilde\psi\Mid F = \psi^F$. If $K/F$ is a finite extension contained in $\wt F$, set $\psi^K = \tilde\psi\Mid K$. 
\enddefinition 
For $\beta\in \Swr{\frak a}$, define 
$$ 
\wt{\scr D}(\frak a,\beta,\tilde\psi) = \underset {K/F} \to {\underset \longleftarrow\to{\text{\rm lim}}}\, \scr D(\frak a_K,\beta,\psi^K), 
\tag 4.2.1 
$$ 
the limit being taken with respect to the canonical restriction maps. The elements of $\wt{\scr D}(\frak a,\beta,\tilde\psi)$ are characters of the group $\r I_\infty(\beta,\frak a)$. Each such character extends uniquely to a continuous character of $\tilde{\r I}(\beta,\frak a)$, so we regard $\wt{\scr D}(\frak a,\beta,\tilde\psi)$ as a set of characters of the group $\tilde{\r I}(\beta,\frak a)$. 
\par 
Let $K/F$ be a finite extension inside $\wt F$. By 2.2 Proposition and the definition (4.2.1), the restriction map 
$$ 
\aligned 
\wt{\scr D}(\frak a,\beta,\tilde\psi) &\longrightarrow \scr D(\frak a_K,\beta,\psi^K), \\ 
\tilde\lambda &\longmapsto \tilde\lambda^K = \tilde\lambda\Mid \r I(\beta,\frak a_K), 
\endaligned 
\tag 4.2.2
$$ 
is surjective. 
\subhead 
4.3 
\endsubhead 
From now on, the following notation will be standard. 
\definition{Notation} 
Set $A = \M nF$ and fix a prime element $\vp_F$ of $F$. Let $\frak a$ be the standard minimal hereditary $\frak o_F$-order in $A$, consisting of all $a\in \M n{\frak o_F}$ that are upper triangular when reduced modulo $\frak p_F$. The order $\frak a$ has a standard prime element $\vP$, such that $\vP^n = \vp_FI_n$. Specifically, all entries $x_{ij}$ of the matrix $\vP$ are zero except $x_{i,i{+}1} = 1$, $1\le i\le n{-}1$, and $x_{n1} = \vp_F$. The Jacobson radical of $\frak a$ is then $\frak p = \vP\frak a = \frak a\vP$. 
\par 
Set $\wt A = A\otimes_F\wt F$ and $\tilde{\frak a} = \frak a\otimes_{\frak o_F}\tilde{\frak o}$. The group $\vO$ acts on $\wt A$ via the first tensor factor, the $F$-algebra $\wt A^\vO$ of $\vO$-fixed points being $A$. Let $\tau$ be a topological generator of the pro-cyclic group $\vO$. Thus $\tau = \sigma_F^z$, for some $z\in \widehat{\Bbb Z}^\times$, where $\widehat{\Bbb Z}$ is the profinite completion of $\Bbb Z$ and $\sigma_F$ is the arithmetic Frobenius. 
\enddefinition 
\proclaim{Proposition} 
Let $m$ be positive divisor of $n$, say $n=md$. The set $B = \wt A^{\tau\vP^m}$ of $\roman{Ad}\,\tau\vP^m$-fixed points in $\wt A$ is a central simple $F$-algebra of dimension $n^2$ and Hasse invariant 
$$ 
\roman{inv}_F\,B = -d^{-1}z+\widehat{\Bbb Z} \in \Bbb Q/\Bbb Z. 
$$ 
In particular, $B \cong \M mD$, for a central $F$-division algebra $D$ of dimension $d^2$. 
\endproclaim 
\demo{Proof} 
Let $\Delta$ be the algebra of diagonal matrices in $A$ and $\wt\Delta$ that in $\wt A$. The vector space $A$ is then the direct sum of the spaces $\Delta\vP^i$, $0\le i < n$, and likewise for $\wt A$. In $\wt A$, each of the spaces $\wt\Delta\vP^i$ is stable under both $\roman{Ad}\,\tau$ and $\roman{Ad}\,\vP$. 
\par 
We deal first with the special case $m = 1$. The $F$-algebra $L = \wt\Delta^{\tau\vP}$ is then an unramified field extension of $F$, of degree $n$. The automorphism $\roman{Ad}\,\vP$ stabilizes $L$, where it acts as the Galois automorphism induced by $\tau^{-1} = \sigma_F^{-z}$. That is, $\sigma_F^{-z}\Mid L = \sigma_F^{-z_0}\Mid L$, for an integer $z_0$ uniquely determined modulo $n$. In other words, $z_0{+}n\widehat{\Bbb Z} = z{+}n\widehat{\Bbb Z}$. The algebra $B$ is thus the classical cyclic division algebra of Hasse invariant $-z_0/n \pmod{\Bbb Z}$: see the Appendix to section 1 in \cite{23}. 
\par 
In the general case, let $e_i \in \Delta$ be the diagonal idempotent matrix with $1$ in the $i$-th place. In particular, $e_i$ is indecomposable and $\bold e = \{e_i:1\le i\le n\}$ is a complete set of orthogonal, indecomposable idempotents in $A$ or in $\wt A$. Viewing the $e_i$ as indexed by the elements of $\Bbb Z/n\Bbb Z$, the automorphism $\roman{Ad}\,\tau$ fixes each $e_i$, while $\roman{Ad}\,\vP$ maps $e_i$ to $e_{i-1}$. 
\par 
Since $\vP^n = \vp_F$, each orbit of $\roman{Ad}\,\vP^m$ on the set $\bold e$ has $d$ elements, where $n = md$, and there are $m$ distinct orbits. For each such orbit $O$, let $e_O$ be the sum of its elements. Thus $e_O \in B = \wt A^{\tau\vP^m}$ and $e_O$ is an idempotent in $B$. The $F$-algebra $e_O B e_O$ has $e_O$ as unit element. Moreover, $\wt F \otimes e_OBe_O = e_O\wt Ae_O$ and $e_O\wt Ae_O \cong \M d{\wt F}$. For, if $\wt V$ is a simple left $\wt A$-module, so that $\wt A = \End{\wt F}{\wt V}$, then $e_O\wt Ae_O = \End{\wt F}{e_O\wt V}$. Consequently, $e_OBe_O$ is a central simple $F$-algebra of dimension $d^2$. 
\par 
On the other hand, the ring $e_O\frak ae_O$ is the standard minimal hereditary order in $e_OAe_O = \M dF$ and its standard prime element is $\vP_O = e_O \vP^me_O$. It satisfies $\vP_O^d = \vp_Fe_O$ and $\ups_F(\det_{e_oBe_o}\vP_O) = 1$. Moreover, $e_OBe_O = (e_O\wt Ae_O)^{\tau\vP_O}$. We are therefore reduced to the special case above and the result follows. \qed 
\enddemo 
\remark{Remarks} 
\roster 
\item  
Since $\vP^n = \vp_F$ and $\tau$ commutes with $\vP$, the operators $\roman{Ad}\,(\tau\vP^m)^n$, $\roman{Ad}\,\tau^n$ are the same. The field $K = \wt F^{\tau^n}$ of $\tau^n$-fixed points in $\wt F$ is of degree $n$ over $F$.  In the notation of the proposition, $\wt A^{\tau^n} = A_K = B_K = \M nK$. Moreover, the set $\frak b_K = \tilde{\frak a}\cap B_K = \frak a_K$ is a minimal hereditary $\frak o_K$-order in $B_K$. 
\item 
We likewise have $\tilde{\frak a}^\tau = \frak a$. As $B = \wt A^{\tau\vP^m}$, so $\frak b = \tilde{\frak a}^{\tau\vP^m} = \frak b_K\cap B$ is a minimal hereditary $\frak o_F$-order in $B$. 
\item 
If $C$ is a central simple $F$-algebra of dimension $n^2$, we may choose the element $\tau$ so that $C \cong \wt A^{\tau\vP^m}$. 
\endroster 
\endremark 
\subhead 
4.4 
\endsubhead 
We write $\wt U = \tilde{\frak a}^\times$ and work in the group $\vO \wt G = \vO\ltimes \wt G$, where $\wt G  = \wt A^\times$.  
\proclaim{Proposition } 
Let $u\in \widetilde U$. There exists $y\in \widetilde U$ such that $u\tau\vP^m = y\tau\vP^my^{-1}$. 
\endproclaim 
\demo{Proof} 
We start with an elementary and familiar observation. 
\proclaim{Cohomological Lemma} 
\roster 
\item 
If $x\in \wt\Bbbk^\times$, there exists $y\in \wt\Bbbk^\times$ such that $x = y^\tau y^{-1}$. 
\item 
If $x\in \wt\Bbbk$, there exists $y\in \wt\Bbbk$ such that $x= y^\tau{-}y$. 
\endroster 
\endproclaim 
\demo{Proof} 
In either part, the element $x$ lies in some finite field $\Bbbk/\Bbbk_F$. In the first part, an elementary argument gives a finite extension $\ell/\Bbbk$ such that $\N\ell{\Bbbk_F}(x) = 1$. If $\vG = \Gal\ell{\Bbbk_F}$, the triviality of the Tate cohomology group $\widehat{\r H}^{-1}(\vG,\ell^\times)$ (Hilbert 90) gives the result. In part (2), we choose $\ell$ so that $\roman{Tr}_{\ell/\Bbbk_F}(x) = 0$. That $\hat{\r H}^{-1}(\vG,\ell) = 1$ implies the result. \qed 
\enddemo 
The group $\tilde{\frak p} = \vP \tilde{\frak a}$ is the Jacobson radical of $\tilde{\frak a}$. Recall that $\tilde{\bs\mu}$ denotes the group of roots of unity in $F_\infty$ of order prime to $p$. Reduction modulo $\vp_F\tilde{\frak o}$ induces an isomorphism $\tilde{\bs\mu} \to \wt\Bbbk^\times$.  Set $\wt U^k = 1{+}\tilde{\frak p}^k$, $k\ge1$. Thus $\wt U$ decomposes as a semi-direct product 
$$ 
\wt U = \tilde{\bs\mu}^n \ltimes \wt U^1. 
$$ 
The groups $\tilde{\bs\mu}^n = \tilde{\bs\mu} \times \tilde{\bs\mu} \times \dots \times \tilde{\bs\mu}$  (with $n$ factors) and $\wt U^1$ are both stable under conjugation by $\tau$ and $\vP^m$ separately.
\proclaim{Lemma 1} 
Let $u\in \tilde{\bs\mu}^n$. There exists $y\in \tilde{\bs\mu}^n$ such that $u\tau \vP^m = y\tau\vP^my^{-1}$. 
\endproclaim 
\demo{Proof} 
The Galois automorphism $\tau$ acts on each factor $\tilde{\bs\mu}$ in the natural way while $\roman{Ad}\,\vP$ permutes them. Write $u = (u_1,\dots ,u_n)$ and likewise for $y$. We have 
$$ 
\tau\vP^m y \vP^{-m}\tau^{-1}m = (\tau(y_{1+m}),\tau(y_{2+m}),\dots,\tau(y_{n+m})), 
$$ 
all subscripts being read modulo $n$. So, if we set $d = n/m$, we have to solve $m$ independent systems of equations in $z_i\in \tilde{\bs\mu}$, of the form 
$$ 
v_i = z_i/\tau(z_{i+1}), \quad 1\le i\le d, 
\tag 4.4.1 
$$ 
for given $v_i\in \tilde{\bs\mu}$. The Cohomological Lemma gives an element $z_1\in \tilde{\bs\mu}$ such that 
$$ 
z_1/\tau^d(z_1) = \prod_{i=1}^{d-1} \tau^{i-1}(v_i). 
$$ 
We solve for $z_j$, $2\le j\le d$, directly from (4.4.1). \qed 
\enddemo 
Given $u\in \wt U$, Lemma 1 shows there exists $y_0\in \wt U$ such that $y_0^{-1}u\tau\vP^my_0 = u_1\tau\vP^m$, for some $u_1\in \wt U^1$. We now set $\wt U^k = 1{+}\tilde{\frak p}^k$, $k\ge1$, and proceed iteratively. 
\proclaim{Lemma 2} 
For an integer $k\ge1$, let $u_k \in \wt U^k$. There exists $y_k \in \wt U^k$ such that 
$$ 
y_k^{-1}u_k\tau\vP^my_k = u_{k+1}\tau\vP^m, 
$$ 
for some $u_{k+1}\in \wt U^{k+1}$. 
\endproclaim 
\demo{Proof} 
Let $\frak q$ be the Jacobson radical of $\frak b$. We have $\tilde{\frak a} = \tilde{\frak o}\frak b$ and $\tilde{\frak p} = \tilde{\frak o}\frak q$. We use the isomorphism 
$$ 
\wt U^k/\wt U^{k+1} \cong \tilde{\frak p}^k/\tilde{\frak p}^{k+1} \cong \frak q^k/\frak q^{k+1} \otimes_{\Bbbk_F} \wt\Bbbk. 
$$ 
The automorphism $\tau\vP^m$ acts trivially on the first tensor factor, and as $\tau$ on the second. Write $u_k = 1{+}x$, where $x\in \tilde{\frak p}^k$ takes the form 
$$ 
x{+}\tilde{\frak p}^{k+1} = \sum_i q_i\otimes \zeta_i, 
$$ 
for a finite number of terms with $q_i\in \frak q^k/\frak q^{k+1}$ and $\zeta_i \in \wt\Bbbk$. The Cohomological Lemma gives $z_i \in \wt\Bbbk$ such that $\zeta_i \equiv z_i-\tau z_i\tau^{-1} \pmod{\tilde{\frak p}^{k+1}}$. Setting 
$$ 
z {+} \wt{\frak  p}^{k+1} = \sum_i q_i\otimes z_i, 
$$ 
the element $y_k = 1{+}z$ has the required property. \qed 
\enddemo 
Using the same notation, set $Y_k = y_0y_1\dots y_k$, $k\ge 0$. The completeness of $\wt F$ ensures that the sequence $\{Y_k\}$ converges to the desired element $y$ of $\wt U$. \qed  
\enddemo 
\subhead 
4.5 
\endsubhead 
Take an algebra $B = \wt A^{\tau\vP^m}$, as in 4.3 Proposition, along with the minimal hereditary $\frak o_F$-order $\frak b = \widetilde{\frak a}^{\tau\vP^m}$ in $B$. Set $H = B^\times$. 
\proclaim{Proposition} 
Let $\beta\in \Swr{\frak a}$, and write $E = F[\beta]$. 
\roster 
\item 
There exists $\gamma\in E$, of valuation $m$, and $y\in \wt U$ satisfying  $y^{-1}E^\times y \i \r K_\frak b \i H$ and $\tau \vP^m = y^{-1}\tau\gamma y$. 
\item 
The element $y$ of \rom{(1)} satisfies 
$$ 
\align 
y^{-1}\tilde{\r I}(\beta,\frak a)y &= \tilde{\r I}(y^{-1}\beta y, \frak b) \quad \text{\rm and} \\ 
\tilde\lambda \circ \roman{Ad}\,y &\in \wt{\scr D}(\frak b,y^{-1}\beta y,\tilde\psi). 
\endalign 
$$ 
\item 
Let $\iota:E\to B$ be an $F$-embedding such that $\iota(E^\times) \i \r K_\frak b$, let $K/F$ be unramified of finite degree divisible by $n$\ and extend $\iota$ to a $K$-isomorphism $\iota_K:A_K \to B_K$. Let $\tilde\lambda\in \wt{\scr D}(\frak a,\beta,\tilde\psi)$ and write 
$$ 
\align 
\tilde\lambda^K &= \tilde\lambda\Mid \r I(\beta,\frak a_K), \\  
(\iota\tilde\lambda^K)^F &= \big(\tilde\lambda^K\circ \iota_K^{-1}\big)\Mid \r I(\iota(\beta),\frak b). 
\endalign 
$$ 
The character $(\iota\tilde\lambda^K)^F$ lies in $\scr D(\frak b,\iota(\beta),\psi^F)$ and 
$$ 
\cind_{\r I(y^{-1}\beta y,\frak b)}^H\big(\tilde\lambda \circ \roman{Ad}\,y\Mid \r I(y^{-1}\beta y,\frak b)\big) \cong \cind_{\r I(\iota(\beta),\frak b)}^H\, \big((\iota\tilde\lambda^K)^F\big). 
$$
\endroster 
\endproclaim  
\demo{Proof} 
Let $\vp$ be a prime element of $E$ and set $\gamma = \vp^m$. There is a unit $u\in \wt U$ such that $\gamma = u\vP^m$. By 4.4 Proposition, there exists $y \in \wt U$ such that $y^{-1}\tau\gamma y = \tau\vP^m$. Since $\vp \in A$, it commutes with $\tau$. Therefore $y^{-1}\vp y$ commutes with $y^{-1}\tau\gamma y = \tau\vP^m$. That is, $y^{-1}\vp y \in H$, whence $y^{-1}Ey \i B$. On the other hand, the group $y^{-1}E^\times y$ normalizes $\tilde{\frak a}$, so it also normalizes $\tilde{\frak a}\cap B = \frak b$. That is, $y^{-1}E^\times y \i \r K_\frak b$. 
\par 
In (2), we have $A_K = B_K$ and $\frak a_K =\frak b_K$. Abbreviating $E_y = y^{-1}Ey$, we have $E_y^\times \i \r K_\frak b$, so $KE_y^\times \i \r K_{\frak b_K} = \r K_{\frak a_K}$. By 3.1 Proposition (2), there exists $x\in \r K_{\frak a_K}$ so that $x^{-1}gx = y^{-1}gy$, for all $g\in KE^\times$. In the language of 3.1, $\iota$ is the embedding $\roman{Ad}\,y^{-1}:E\to B$ and $\roman{Ad}\,x^{-1}$ is the extension $\iota_K$ of $\iota$ to a $K$-isomorphism $A_K \to B_K = A_K$. The first assertion in (2) now follows from 3.1 Proposition on passing to the limit over $K$. The second assertion of (2) follows the same course. 
\par 
In (3), we use 3.1 Proposition again: $\iota$ extends to a $K$-automorphism $\iota_K$ of $A_K$, stabilizing $\frak a_K$. This has the form $\iota_K = \roman{Ad}\,y_0$, for some $y_0 \in U_{\frak a_K}$. So, by definition, $\pi_H((\iota\tilde\lambda^K)^F)$ is induced by $\tilde\lambda\circ\roman{Ad}\,y_0 \Mid \tilde{\r I}(y_0^{-1}\beta y_0, \frak b)\cap H$. But, as in 3.1 Proposition, $\roman{Ad}\,y = \roman{Ad}\, zy_0x$, for some $x\in H$ and $x\in \wt F[\beta]^\times$. The factor $\roman{Ad}\,z$ has no effect on the inducing datum, while $\roman{Ad}\,x$ does not change the equivalence class of the induced representation. \qed 
\enddemo 
\head 
5. Basic character relation 
\endhead 
We prove our pivotal result. We use the notation of 4.3 along with a central simple $F$-algebra $B$ of dimension $n^2$, realized in the form $B = \wt A^{\tau\vP^m}$ as in 4.5. We set $H = B^\times$ and $\frak b = \tilde{\frak a}\cap B$, again as in 4.5. Recall that, in this scheme, $A = \M nF$ and $G = \GL nF$.   We assume throughout that $p\neq 2$. 
\subhead 
5.1 
\endsubhead 
We evaluate characters of representations of $G$ and $H$ at a certain class of elements as follows. 
\definition{Definition} 
Let $H^\roman{wr}$ be the set of elements $h$ of $H$ satisfying the following conditions: 
\roster 
\item 
$\ups_F(\det_B(h))$ is not divisible by $p$ and 
\item 
there exists a minimal hereditary $\frak o_F$-order $\frak b_1$ in $B$ such that $h\in \r K_{\frak b_1}$. 
\endroster 
\enddefinition 
Let $h\in H^\roman{wr}$. The algebra $L = F[h]$ is a field, totally ramified of degree $n$ over $F$, and $L^\times \i \r K_{\frak b_1}$. In particular, the reduced characteristic polynomial $\roman{ch}_B(t;h) \in F[t]$ of $h$ is irreducible over $F$. Let $H^\roman{wr}_\roman{reg}$ be the set of $h\in H^\roman{wr}$ for which $\roman{ch}_B(t;h)$ is also separable. Thus an element of $H^\roman{wr}$ is ``elliptic quasi-regular'', in the sense of \cite{4} A.2, while any $h\in H^\roman{wr}_\roman{reg}$ is elliptic regular in the customary sense. The sets $H^\roman{wr}$, $H^\roman{wr}_\roman{reg}$ are stable under conjugation by $H$
\par 
The sets $G^\roman{wr}$, $G^\roman{wr}_\roman{reg}$ are defined in the same way. 
\proclaim{Lemma} 
Let $g\in G^\roman{wr}$. There is a unique $H$-conjugacy class of elements $h\in H^\roman{wr}$ such that $\roman{ch}_B(t;h) = \roman{ch}_A(t;g)$. Equality of reduced characteristic polynomials induces canonical bijections 
$$ 
\aligned 
\roman{Ad}\,H\backslash H^\roman{wr} &\longrightarrow \roman{Ad}\,G\backslash G^\roman{wr}, \\ 
\roman{Ad}\,H\backslash H^\roman{wr}_\roman{reg} &\longrightarrow \roman{Ad}\,G\backslash G^\roman{wr}_\roman{reg}. 
\endaligned 
\tag 5.1.1 
$$ 
\endproclaim 
The proof is immediate. The lemma clearly remains valid on replacing $G = \GL nF$ by an inner form. We refer to the bijections (5.1.1), and their inverses, as {\it association.} 
\subhead 
5.2 
\endsubhead  
We use the following additional notation. 
\definition{Notation} 
\roster 
\item 
Let $\beta\in \Swr{\frak a}$, let $\lambda\in \scr D(\frak a,\beta,\psi^F)$ and write $\pi_G = \pi_G(\lambda)$ for the representation of $G$ induced by $\lambda$, as in 2.3. 
\item 
Let $\tilde\lambda \in \wt{\scr D}(\frak a,\beta,\tilde\psi)$ satisfy $\tilde\lambda\Mid \r I(\beta,\frak a) = \lambda$. Let $y\in \wt U$ satisfy 4.5 Proposition. Let $\pi_H$ be the representation of $H$ induced by the character $\tilde\lambda\circ\roman{Ad}\,y\Mid \r I(y^{-1}\beta y,\frak b)$. 
\endroster 
\enddefinition 
\proclaim{Theorem} 
If $g\in G^\roman{wr}$ is associate to $h\in H^\roman{wr}$, then 
$$ 
\roman{tr}\,\pi_G(g) = \roman{tr}\,\pi_H(h). 
\tag 5.2.1 
$$ 
\endproclaim 
The proof occupies the rest of the section. 
\subhead 
5.3 
\endsubhead 
We work first with the group $H$. With notation as in 5.2, write 
$$ 
\align 
\tilde\kappa &= \tilde\lambda \circ \roman{Ad}\, y \in \wt{\scr D}(\frak b,y^{-1}\beta y, \tilde\psi), \\ 
\kappa &= \tilde\kappa\Mid \r I(y^{-1}\beta y,\frak b).
\endalign 
$$ 
Thus $\kappa \in \scr D(\frak b,y^{-1}\beta y,\psi^F)$ and $\pi_H$ is induced by $\kappa$. Abbreviate $\r I(y^{-1}\beta y,\frak b) =I_H$ and $\wt{\r I}(y^{-1}\beta y,\frak b) = \tilde I_H$. 
\par 
Let $h\in H^\roman{wr}$. The Mackey induction formula gives 
$$ 
\roman{tr}\,\pi_H(h) = \sum\Sb x\in H/I_H, \\ x^{-1}hx\in I_H \endSb \kappa(x^{-1}hx). 
$$ 
The condition $x^{-1}hx\in I_H$ is equivalent to $hxI_H = xI_H$, that is, to $xI_H$ being a fixed point for the natural left translation action of $h$ on $H/I_H$. We may therefore re-write this character expansion in the form 
$$ 
\roman{tr}\,\pi_H(h) = \sum_{x\in (H/I_H)^h} \kappa(x^{-1}hx) .
\tag 5.3.1 
$$ 
\remark{Remark} 
Since $h$ is elliptic quasi-regular, the argument of \cite{6} 1.2 Lemma applies to show that the expansion (5.3.1) has only finitely many terms. Consequently, there are no convergence issues to be considered. 
\endremark 
The natural conjugation action of $\tau\vP^m$ on $\wt G$ stabilizes $\tilde I_H$, so $\tau\vP^m$ acts on the coset space $\wt G/\tilde I_H$ by conjugation. We write $\big(\wt G/\tilde I_H\big)^{\tau\vP^m}$ for the set of fixed points. On the other hand, $\tilde I_H^{\tau\vP^m} = \tilde I_H\cap H = I_H$. 
\proclaim{Lemma} 
The canonical map $H/I_H \to \big(\wt G/\tilde I_H\big)^{\tau\vP^m}$ is a bijection. 
\endproclaim 
\demo{Proof} 
The map in question is surely injective. For $x\in \wt G$, the coset $x\tilde I_H$ is $\tau\vP^m$-fixed if and only if $\tau\vP^mx = xj\tau\vP^m$, for some $j\in \tilde I_H$. Any such element $j$ must lie in $\tilde I^0_H = \wt U\cap \tilde I_H$. There exists $k\in \tilde I^0_H$ such that $j\tau\vP^m = k\tau \vP^m k^{-1}$: this is proved is the same way as 4.4 Proposition but is easier since $\tilde I^0_H = \tilde{\bs\mu}\times \tilde I^1_H$. For this element $k$, we have 
$$ 
\tau\vP^mx = xk\tau \vP^mk^{-1}, \quad\text{or}\quad \tau\vP^mxk =xk\tau\vP^m, 
$$ 
giving $xk\in H$, as desired. \qed 
\enddemo 
The obvious inclusion of $\wt G$ in $\langle \tau\vP^m \rangle \ltimes \wt G$ induces a bijection 
$$ 
\wt G/\tilde I_H = \langle \tau\vP^m \rangle \ltimes \wt G\big/ \langle \tau\vP^m \rangle \ltimes \tilde I_H. 
$$ 
We use it to extend the translation action of $\wt G$ on $\wt G/\tilde I_H$ to one of $\langle \tau\vP^m \rangle \ltimes \wt G$. The set of $\tau\vP^m$-fixed points in $\wt G/\tilde I_H$ for this action is $H/I_H$ (by the lemma), so we may re-write (5.3.1) as 
$$ 
\roman{tr}\,\pi_H(h) = \sum_{x\in (\wt G/\tilde I_H)^{\langle h, \tau\vP^m \rangle}} \tilde\kappa(x^{-1}h x), \quad h\in H^\roman{wr}. 
\tag 5.3.2 
$$ 
\subhead 
5.4 
\endsubhead 
We apply the same argument to the representation $\pi_G$. Abbreviating $\tilde I = \tilde{\r I}(\beta,\frak a)$, we find 
$$ 
\roman{tr}\,\pi_G(g) = \sum_{x\in (\wt G/\tilde I)^{\langle g,\tau\rangle}} \tilde\lambda(x^{-1}gx), \quad g\in G^\roman{wr}. 
\tag 5.4.1 
$$ 
When evaluating this finite sum, there is no loss entailed in assuming $g\in \r K_\frak a$. 
\proclaim{Lemma 1} 
Let $g\in \r K_\frak a\cap G^\roman{wr}$and let $g_0 \in F[g]$ have valuation $\ups_{F[g]}(g_0) = m$. There exists $t \in \wt U$ such that $t^{-1} F[g]^\times t \i \r K_\frak b$ and $t^{-1}g_0\tau t = \tau\vP^m$. 
\endproclaim 
The proof is identical to that of 4.5 Proposition (1), so we say no more of it. The element $h = t^{-1}gt$ lies in $\r K_\frak b\cap H^\roman{wr}$ and is associate to $g$. We therefore compare (5.4.1) with (5.3.2) evaluated at this element $h$. 
\par 
We relate the index sets $(\wt G/\tilde I)^{\langle\tau,g\rangle}$, $(\wt G/\tilde I_H)^{\langle \tau\vP^m, h \rangle}$. To do this, we view $\vO\ltimes \wt G$ as acting, by left translation, on $\wt G/\tilde I$ and $\langle \tau\vP^m\rangle \ltimes \wt G$ likewise on $\wt G/\tilde I_H$. 
\proclaim{Lemma 2} 
The map 
$$ 
\align 
\vF: \wt G/\tilde I &\longrightarrow \wt G/\tilde I_H, \\ 
x\tilde I &\longmapsto t^{-1}xy \tilde I_H, 
\endalign 
$$ 
induces a bijection 
$$ 
\big(\wt G/\tilde I\big)^{\langle\tau,g\rangle} @>{\ \ \approx\ \ }>> \big(\wt G/\tilde I_H\big)^{\langle \tau\vP^m, h \rangle}. 
$$ 
\endproclaim 
\demo{Proof} 
The defining properties of $y\in \wt U$ ({\it cf\.} 5.2 Notation) are 
$$ 
y^{-1}E^\times y \i \r K_\frak b \quad\text{and} \quad \tau\vP^m = y^{-1}\tau\gamma y, 
$$ 
for a certain element $\gamma$ of $E$ as in 4.5 Proposition. We also have $y^{-1}E^\times y \i \r I(y^{-1}\beta y,\frak b) = I_H \i \tilde I_H$. 
\par 
Likewise, $t \in \wt U$ satisfies 
$$ 
t^{-1}F[g]^\times t \i \r K_\frak b \quad \text{and} \quad t^{-1}g_0\tau t  = \tau\vP^m. 
$$ 
We have set $h = t^{-1}gt$. 
\par 
Immediately, the map $\vF$ is a bijection $\wt G/\tilde I \to \wt G/\tilde I_H$. Let $x\in \wt G$ and suppose that $x\tilde I$ is fixed by $\tau$ and $g$. By the obvious analogue of 5.3 Lemma, we may assume $x\in G$. Since $gxI = xI$, the element $x$ conjugates $g$ into $I$. The algebra $F[g]$ is a field and $g$ in minimal over $F$, in the sense of \cite{9} (1.4.14), whence it follows readily that $x$ conjugates $F[g]^\times$ into $I$. Thus 
$$ 
h\vF(x\tilde I) = h t^{-1}xy\tilde I_H = t^{-1}gxy \tilde I_H = \vF(gx\tilde I) = \vF(x\tilde I), 
$$ 
as desired. 
\par 
Now consider 
$$ 
\align 
\tau\vP^m\vF(x\tilde I) &= \tau\vP^m t^{-1}xy\tilde I_H = t^{-1}g_0\tau xy\tilde I_H \\ 
&= t^{-1}g_0x\tau y\tilde I_H = t^{-1}g_0x\gamma^{-1}y\tilde I_H\tau\vP^m. 
\endalign 
$$ 
Since $\gamma \in E$, we have $y^{-1}\gamma^{-1} \in I_H$, whence 
$$ 
\tau\vP^m\vF(x\tilde I) = t^{-1}g_0xy\tilde I_H\tau\vP^m = \vF(g_0x\tilde I)\tau\vP^m = \vF(x\tilde I)\tau\vP^m.  
$$ 
Thus $\vF(x\tilde I)$ is fixed by $\tau\vP^m$ and $h$, as required. The argument is reversible and the lemma is proven. \qed 
\enddemo 
\subhead 
5.5 
\endsubhead 
We prove 5.2 Theorem. Let $xI \in (G/I)^g$. The bijection $\vF$ of 5.4 Lemma 2 gives a coset $x'I_H = t^{-1}xyj(x) I_H \in (H/I_H)^h$, for some $j(x) \in \tilde I_H$ uniquely determined modulo $I_H$. The contribution to (5.3.1) from the coset $x'I_H$ is 
$$ 
\align 
\kappa({x'}^{-1}hx') &= \tilde\kappa({x'}^{-1}hx') \\ 
&=  \tilde\kappa(j(x)^{-1}y^{-1}x^{-1}tht^{-1}xyj(x)). 
\endalign 
$$ 
Since $\tilde\kappa$ is a character of $\tilde I_H$, it is invariant under conjugation by $j(x)$. Recalling that $\tilde\kappa = \tilde \lambda \circ \roman{Ad}\,y$, this expression reduces to 
$$ 
\kappa({x'}^{-1}hx') = \tilde\lambda(x^{-1}tht^{-1}x) = \lambda(x^{-1}gx). 
$$ 
Lemma 2 now implies $\roman{tr}\,\pi_G(g) = \roman{tr}\,\pi_H(h)$, as required. \qed 
\head 
6. Consequences 
\endhead 
We derive from 5.2 Theorem the main results of the paper. 
\subhead 
6.1 
\endsubhead 
Let $G$, $H$ be inner forms of $\GL nF$, and let 
$$ 
T_G^H: \awr G \longrightarrow \awr H 
$$ 
denote the bijection induced by the Jacquet-Langlands correspondence. 
\proclaim{Theorem} 
Let $G$ and $H$ be inner forms of $\GL nF$. Let $\pi\in \awr G$ and let $\pi'\in \awr H$ be a parametric transfer of $\pi$. 
\roster 
\item 
The representations $\pi$, $\pi'$ are related by 
$$ 
\pi' = T_G^H(\pi). 
$$ 
\item 
If $\rho\in \awr H$ satisfies 
$$ 
\roman{tr}\,\rho(h) = \roman{tr}\,\pi(g), 
$$ 
for all $g\in G^\roman{wr}_\roman{reg}$ with associate $h\in H^\roman{wr}_\roman{reg}$, then $\rho = T_G^H(\pi) = \pi'$. 
\endroster 
\endproclaim 
Consequently, the equivalence class of a parametric transfer $\pi'$ of $\pi$ depends only on that of $\pi$, and not on the choices made in the definition in 3.4. The theorem implies the endo-class transfer theorem in this special case. 
\proclaim{Corollary} 
Let $\pi\in \awr G$ and let $\rho  = T_G^H(\pi)\in \awr H$. If $\theta_\pi$, $\theta_\rho$ are simple characters contained in $\pi$, $\rho$ respectively, then $\theta_\pi$, $\theta_\rho$ are endo-equivalent. 
\endproclaim 
The proofs are in 6.4, following some preparatory material in 6.3. 
\subhead 
6.2 
\endsubhead 
Before proceeding to the proofs, we return to the discussions leading to 3.5 in order to write the Jacquet-Langlands correspondence in explicit form. Re-setting the notation, let $A$ and $B$ be central simple $F$-algebra of dimension $n^2$, let $\frak a$ and $\frak b$ be minimal hereditary $\frak o_F$-orders in $A$ and $B$ respectively and write $G = A^\times$, $H = B^\times$. 
\par 
Let $\beta\in \Swr{\frak a}$, write $E = F[\beta]$ and let $\theta\in \aC(\frak a,\beta,\psi^F)$. Let $\xi$ be a character of $E^\times$ such that $\xi\Mid U^1_E = \theta\Mid U^1_E$. Let $\iota:E\to B$ be an $F$-embedding such that $\iota(E^\times)\i \r K_\frak b$, extended to an isomorphism $\iota_K:A_K \to B_K$. Define $\iota\theta^F_K \in \aC(\frak b,\iota(\beta),\psi^F)$ as in (3.3.1).  
\proclaim{Theorem} 
Abbreviating $\tau = \iota\theta^F_K \in \aC(\frak b,\iota(\beta),\psi^F)$, we have 
$$ 
T_G^H\,\pi_G(\xi\odot \theta) = \pi_H((\xi\circ\iota^{-1})\odot\tau). 
\tag 6.2.1 
$$ 
\endproclaim 
\demo{Proof} 
The definitions ensure that $\pi_H((\xi\circ\iota^{-1})\odot\tau)$ is a parametric transfer of $\pi_G(\xi\odot \theta)$ (3.5 Proposition), so the result follows from 6.1 Theorem. \qed 
\enddemo 
\subhead 
6.3 
\endsubhead 
In this sub-section, the algebra $B$ will be {\it a division algebra.} We need some special properties of the characters of irreducible smooth representations of $H = B^\times$. 
\proclaim{Lemma} 
Let $\pi\in \awr H$ and let $\rho$ be an irreducible smooth representation of $H$. The following conditions are equivalent. 
\roster 
\item 
$\roman{tr}\,\pi(g) = \roman{tr}\,\rho(g)$ for all $g\in H^\roman{wr}_\roman{reg}$. 
\item 
$\roman{tr}\,\pi(g) = \roman{tr}\,\rho(g)$ for all $g\in H^\roman{wr}$. 
\item 
$\rho\cong \pi$. 
\endroster 
\endproclaim 
\demo{Proof} 
The characters $\roman{tr}\,\pi$, $\roman{tr}\,\rho$ are locally constant functions on $H$. The set $H^\roman{wr}_\roman{reg}$ is dense in $H^\roman{wr}$. Thus (1) is equivalent to (2) and surely (3) implies (1). 
\par 
Let $X$ be the group of unramified characters $\chi$ of $F^\times$ such that $\chi^n=1$. Let $h\in H^\roman{wr}$, $\ups_F(\det_B(h)) = 1$. We form the function 
$$ 
\align 
\Phi_\pi(g) &= n^{-1}\sum_{\chi\in X} \chi({\det}_B(gh^{-1})\,\roman{tr}\,\pi(g) \\ 
&= n^{-1}\sum_{\chi\in X} \chi({\det}_B(h^{-1}))\, \roman{tr}\,\chi\pi(g), \quad g\in H. 
\endalign 
$$ 
Define $\Phi_\rho$ similarly. Both functions $\Phi_\pi$, $\Phi_\rho$ are supported in the set of $g\in H$ with $\ups_F(\det_B(g)) \equiv 1 \pmod n$. This set is contained in $H^\roman{wr}$ so (2) implies $\Phi_\pi(g) = \Phi_\rho(g)$ for all $g\in H$. The set of characters of irreducible smooth representations of $H$, viewed as functions on $H$, is linearly independent. We conclude that $\rho = \chi\pi$, for some $\chi\in X$ with $\chi(\det_B(h)) = 1$. That is, $\chi = 1$ so $\rho\cong\pi$, as required for (3). \qed 
\enddemo 
\subhead 
6.4 
\endsubhead 
We prove the theorem. Since the Jacquet-Langlands correspondence and parametric transfer are transitive, it is enough to prove the theorem and the corollary under the assumption $G = \GL nF$. 
\par 
Initially take $H = B^\times$, for a division algebra $B$. Let $\pi \in \awr G$ and let $\rho \in \awr H$ be a parametric transfer of $\pi$. Set $\rho' = T_G^H(\pi)$. Let $g \in G^\roman{wr}_\roman{reg}$ and let $h \in H^\roman{wr}_\roman{reg}$ be associate to $g$. Since $n = p^r$, $p\neq 2$, we have $\roman{tr}\,\rho'(h) = \roman{tr}\,\pi(g)$ by definition. However, 5.2 Theorem and 4.5 Proposition (3) together  give $\roman{tr}\,\pi(g) = \roman{tr}\,\rho(h)$ so 6.3 Lemma implies $\rho = \rho'$. This yields an intermediate conclusion : 
\proclaim{Lemma} 
Let $\pi_1, \pi_2 \in \awr G$. If $\roman{tr}\,\pi_1(g) = \roman{tr}\,\pi_2(g)$ for all $g \in G^\roman{wr}_\roman{reg}$, then $\pi_1 = \pi_2$. 
\endproclaim 
We pass to the case where $H$ is arbitrary. Let $\pi\in \awr G$, let $\rho\in \awr H$ be a parametric transfer of $\pi$ and let $\pi'$ be the unique element of $\awr G$ for which $\rho = T_G^H\,\pi'$. Let $g \in G^\roman{wr}_\roman{reg}$ and let $h \in H^\roman{wr}_\roman{reg}$ be associate to $g$. By 5.2 Theorem and 4.5 Proposition (3) again, 
$$ 
\roman{tr}\,\pi(g) = \roman{tr}\,\rho(h) = \roman{tr}\,\pi'(g). 
$$ 
The relation $\roman{tr}\,\pi(g) = \roman{tr}\,\pi'(g)$ holds for all $g\in G^\roman{wr}_\roman{reg}$, and the lemma implies $\pi = \pi'$. This completes the proof of the theorem. \qed 
\par 
The corollary now follows from 3.3 Proposition. \qed 
\remark{Remark} 
We could have argued here in terms of elliptic quasi-regular elements $g\in G^\roman{wr}$. However, elliptic regular elements suffice to give the result and the extra precision can be useful in the context of linear independence of characters. 
\endremark 
\remark{\bf Correction} 
On this subject, 3.1 Corollary 3 of \cite{5} is wrong for rather trivial reasons. For a counterexample, take $\pi \in \awr G$ and consider the set of representations $\chi\pi$, as $\chi$ ranges over all unramified characters of $F^\times$ such that $\chi^n = 1$. The set of characters $\roman{tr}\,\chi\pi$ is then linearly dependent on both $G^\roman{wr}_\roman{reg}$ and $G^\roman{wr}$. This is essentially the only counterexample. The error has no effect on either \cite{5} or this paper.  
\endremark 
 
\enddocument